\documentclass{article}
\usepackage{graphicx} 
\usepackage{amsmath}
\usepackage{mathtools}
\usepackage{algorithm}
\usepackage{algpseudocode}
\usepackage[square,numbers,sort&compress]{natbib}
\usepackage{xcolor}
\usepackage{hyperref}
\usepackage{makecell}
\usepackage{subcaption}
\usepackage{authblk}

\title{Time-vectorized numerical integration for systems of ODEs}
\author[1]{Mark C. Messner\footnote{Corresponding author, email \url{messner@anl.gov}}}
\author[1]{Tianchen Hu}
\author[1,2]{Tianju Chen}
\affil[1]{Argonne National Laboratory, Applied Materials Division}
\affil[2]{Presently at the University of Shanghai for Science and Technology, School of Materials and Chemistry}
\date{September 2023}

\begin{document}

\maketitle

\begin{abstract}
   Stiff systems of ordinary differential equations (ODEs) and sparse training data are common in scientific problems.  This paper describes efficient, implicit, vectorized methods for integrating stiff systems of ordinary differential equations through time and calculating parameter gradients with the adjoint method.  The main innovation is to vectorize the problem both over the number of independent times series and over a batch or ``chunk'' of sequential time steps, effectively vectorizing the assembly of the implicit system of ODEs. The block-bidiagonal structure of the linearized implicit system for the backward Euler method allows for further vectorization using parallel cyclic reduction (PCR).  Vectorizing over both axes of the input data provides a higher bandwidth of calculations to the computing device, allowing even problems with comparatively sparse data to fully utilize modern GPUs and achieving speed ups of greater than 100x, compared to standard, sequential time integration.  We demonstrate the advantages of implicit, vectorized time integration with several example problems, drawn from both analytical stiff and non-stiff ODE models as well as neural ODE models.  We also describe and provide a freely available open-source implementation of the methods developed here.
\end{abstract}

\section{Introduction}

\subsection{Review}

Differential equation models are a fundamental scientific tool for explaining physical phenomenon.  Developing new ordinary differential equation (ODE) models to represent physical observations is a regular part of scientific progress. Formulating ODE models generally has two steps: deriving, often based on physical principles, the form of the ODE model and then training the model parameters against experimental observations.  More recently, neural ODEs have emerged as a data-driven approach, reducing the need for physics-based derivations of new model forms.  Here, a neural network (NN) replaces an analytic model, which can then be trained with the experimental data to determine the model form by finding the optimal parameters.

The paper sparking the recent interest in neural ODEs \cite{chen2018neural} envisioned a broad set of applications, but much recent work focuses on using neural ODEs to model physical phenomena (e.g. \cite{dandekar2020bayesian, brucker2021grey, puliyanda2023benchmarking, lai2021structural, roehrl2020modeling,portwood2019turbulence} among many others).  Researchers have identified several challenges in applying neural ODEs to physical behavior \cite{kidger2022neural}, including difficulty in finding model forms that converge to match physical conservation laws \cite{pmlr-v162-zhu22f,rehman2023physics,hansen2023learning}, uncertainty in formal analysis and difficulties in numerical convergence for neural ODEs \cite{pmlr-v162-zhu22f,gusak2020towards,ott2020resnet}, and challenges with neural ODEs targeting stiff problems \cite{kim2021stiff,de2022physics}.  Stiff ODEs are defined, somewhat nebulously, as systems of ODEs that are difficult to integrate with explicit numerical time integration methods. 

This paper addresses two key issues for efficiently training ODE models, including neural ODEs, on GPUs and other accelerators.  Both of the challenges addressed here are particularly stark when training models for physical problems, but may apply to more general neural ODE models as well.  The two key issues are:
\begin{enumerate}
    \item Current vectorized time integration methods focus on explicit numerical integration schemes.  We are not aware of any numerically efficient implementations of \emph{implicit} time integration methods for solving inherently stiff problems.
    \item Training approaches relying on ODE integration do not perform well for data-sparse problems where vectorizing only over a single batch dimension (i.e. the number of time series) does not provide sufficient bandwidth to efficiently use modern GPUs.
\end{enumerate}

Stiff problems are fairly common in scientific applications and there is a need for efficient, implicit integration methods for both classical analytical and neural models.  Researchers may want to apply  machine learning optimization techniques to classical models \cite{chen2023training} and past work suggests there is a set of inherently-stiff phenomena that seem to require implicit methods to accurately model with neural ODEs \cite{kim2021stiff}.  Extending the basic approach in \citet{chen2018neural} to implicit methods is trivial --- several past works mention implementations in passing --- but \emph{efficiently} implementing vectorized implicit methods is somewhat more challenging.  

The second issue is perhaps more interesting.  Many problems that could be addressed with neural ODEs are data-sparse.  Here a data-sparse problem is one where vectorizing over the number of experimental observations (time series or points from a time series) does not provide sufficient parallel work to fully subscribe the available computational resources.  For physically motivated time series problems we often only have 10s or 100s of time series available for training.  Vectorizing only over the number of independent time series does not fully utilize even a single modern GPU.

\begin{figure}[!htb]
    \centering
    \includegraphics[width=\textwidth]{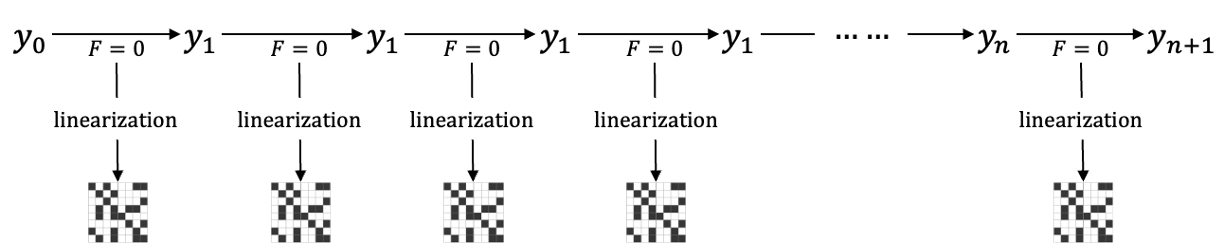}
    \caption{Schematic illustration of the general sequential time integration algorithm for an implicit method. The state is numerically integrated from $y_0$ to $y_{n+1}$. The dynamics is governed by the ODE $F=0$ (see Section~\ref{sec:algs} for a formal definition). The implicit time integration method relies on the linearizations of the system of ODEs for each time increment, which are illustrated as grayscale sparsity patterns. }
    \label{fig:serial_ti}
\end{figure}

\begin{figure}[!htb]
    \centering
    \includegraphics[width=\textwidth]{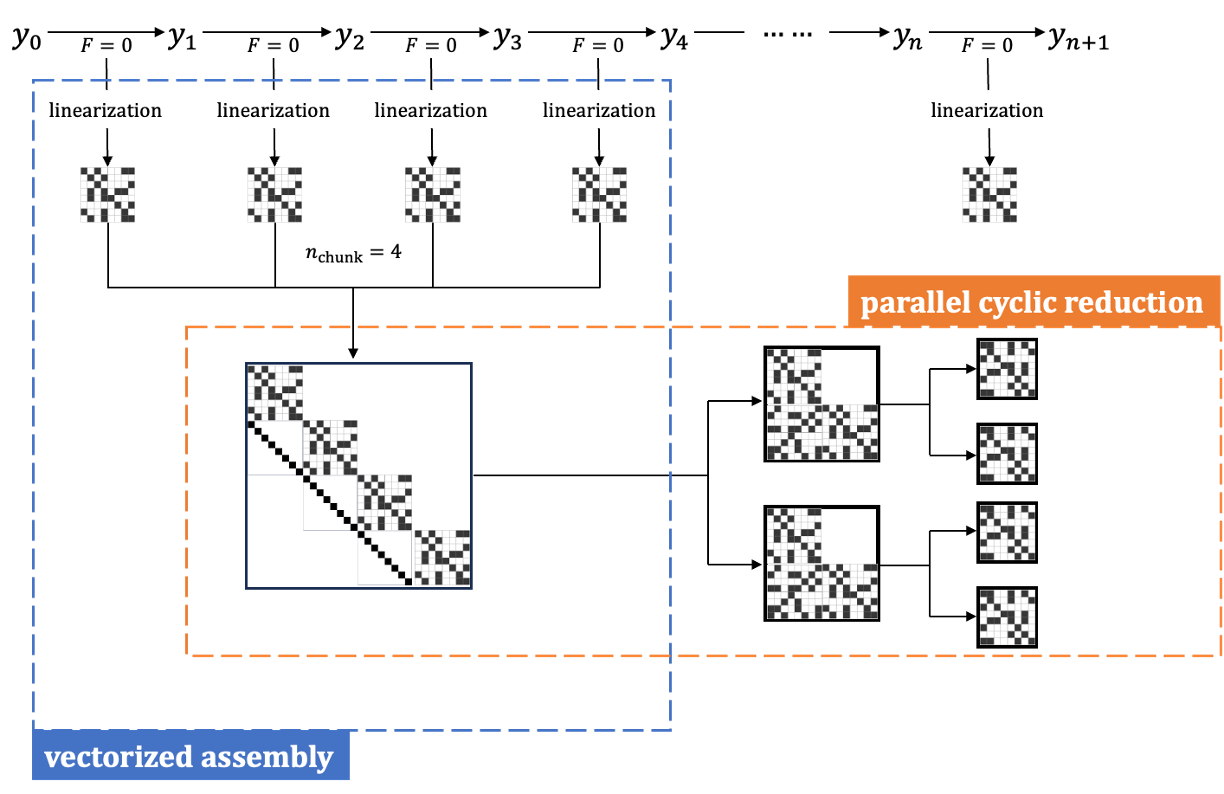}
    \caption{Schematic illustration of the proposed vectorized time integration algorithm, again for an implicit method. In contrast to the sequential time integration algorithm illustrated in Figure~\ref{fig:serial_ti}, the proposed algorithm not only vectorizes the assembly of the linear system of equations, but can also vectorize the matrix factorization using parallel cyclic reduction (Section~\ref{sec:pcr}) by exploiting the special block-bidiagonal structure. }
    \label{fig:vectorized_ti}
\end{figure}

The solution is to vectorize another dimension of the problem --- the number of time steps in the time series to integrate --- effectively ``parallelizing'' time integration.  A first glance this would seem to be impossible, as numerical time integration is an inherently sequential algorithm. The general flow of time integration algorithms is illustrated in Figure~\ref{fig:serial_ti}.  However, inspired by time-space finite element approaches \cite{hughes1988space}, we can increase the computational bandwidth presented to the GPU by integrating a batch or ``chunk''\footnote{To try to mitigate confusion, we refer to a ``batch'' of independent times series and a ``chunk'' of time steps to integrate at once within each time series.} of time steps in a single pass.  This approach not only presents better opportunities to vectorize the computational work associated with evaluating the model and, for implicit methods, the model Jacobian, but also can provide true parallelization opportunities through the special structure of the resulting batched, blocked bidiagonal linear operators arising from chunking the forward and backward Euler integration methods.  A graphical summary of the proposed algorithm is illustrated in Figure~\ref{fig:vectorized_ti}.  The optimal number of time steps to chunk depends on the problem and the computing device, but we demonstrate speed ups in excess of 100x, comparing sequential to chunked time integration, for representative sample problems.

\subsection{Organization}

Section \ref{sec:algs} describes the vectorized time integration algorithm.  We focus on the implicit backward Euler method, though the approach also works for explicit methods and other numerical integration techniques beyond the Euler methods.  This section focuses on the forward pass required to integrate the model through time, an appendix provides the details of the adjoint backward pass.  Vectorized time integration for implicit methods, using Newton's method to solve the resulting discrete nonlinear equations, requires solving a large linear system with a special matrix structure.  Section \ref{sec:bidiagonal} describes two potential algorithms for efficiently solving these linear systems on GPUs.  Section \ref{sec:implementation} includes a few notes on our implementation of the entire approach.  We provide an open-source implementation available in the \emph{pyoptmat} python package.\footnote{https://github.com/Argonne-National-Laboratory/pyoptmat}.  Solving the discrete algebraic equations resulting from implicit methods with Newton's method requires the Jacobian of the ODE.  For some models, i.e. neural ODEs, obtaining this Jacobian analytically can be challenging.  This section also describes the implementation with vectorized automatic differentiation (AD) approaches to calculate the Jacobian.

Section \ref{sec:example-sec} describes several examples demonstrating the computational efficiency of the algorithms and their implementation in \emph{pyoptmat}.  These examples include both non-stiff and stiff physically motivated problems and pure neural ODEs.  The examples also illustrate the cost of calculating the Jacobian for implicit methods with AD, compared to an analytic implementation.  Finally, Section \ref{sec:conclusions} summarizes the key results presented here.

\subsection{Our contributions}
\begin{enumerate}
    \item A numerical time integration strategy that vectorizes over both the number of time series and some chunk of time integration steps.  The method can be applied to implicit or explicit integration schemes, but are most effective for implicit time integration. The method provides speed ups, counting both the forward and backward pass, of greater than 100x for data-sparse problems compared to sequential algorithms batched only over the number of independent time series.
    \item An efficient GPU implementation of implicit time integration schemes with sufficient performance to train ODE and neural-ODE models with conventional, gradient-descent training algorithms.
\end{enumerate}

\section{Chunked time integration \label{sec:algs}}

\subsection{Basic definitions for batched systems of ODEs}

This work considers batched time integration and gradient calculations for systems of ordinary differential equations.  For a system of ODEs
\begin{alignat}{3}
    \dot{y} & = & h\left(y, t ; p \right) \\
    y\left(0\right) & = & y_0
\end{alignat}
with $\left|y\right| = n_{size}$ the problem size, $t$ time, and $p$ the problem parameters define time integration as the problem of finding 
\begin{equation}
    y_i = \int_{0}^{t_i} \dot{y} dt
    \label{eq:time-integration}
\end{equation}
for $n_{time}$ points in time $\left[t_1, t_2, \dots, t_{n_{time}}\right]$ and the gradient calculation as the problem of finding the derivative of a function of one of these results $L(y_i)$ with $L$ some arbitrary scalar-valued function with respect to the problem parameters
\begin{equation}
    g_i = \frac{d L(y_i)}{d p}.
    \label{eq:gradient}
\end{equation}

Almost always, time integration algorithms must proceed incrementally through a series of time points rather than integrate directly to a time of interest, as in Eq. \ref{eq:time-integration}.  For ODEs describing physical systems, users are often interested in examining the solution trajectory at intermediate time points.  Moreover, even if only the final values are of interest, numerical integration schemes have finite accuracy, meaning obtaining an accurate solution often requires marching incrementally through smaller time steps to reach the desired solution.

Batched time integration and gradient calculations mean that we might consider a collection of $n_{batch}$ sets of 
initial conditions $y_{0}^{(j)}$, parameters $p^{(j)}$, forcing functions $f^{(j)}$, and even time series $t_{i}^{(j)}$.  The goal for a batched time integration and gradient calculation algorithm is to simultaneously calculate the solution to these $n_{batch}$ ODEs in a vectorized manner.  While the initial conditions, parameters, forcing functions, and time points all may vary across the batch, the basic model function $h$ remains the same.  As indicated here, subscripts indicate different points in the time series of solutions,  parenthetical superscripts indicate different solutions in the batch, and we leave the problem size $n_{size}$ implicit in $\left|y\right|$.  The manuscript omits the superscripts and subscripts where they are not required for clarity.

\subsection{Implicit time integration schemes}

\subsubsection{Time integration}

\paragraph{Implicit methods}

This paper considers implicit numerical time integration schemes, where the solution of the ODE at the next point in time $t_{i+1}$ is given as the solution to an implicit function involving the previous solution and the next solution:
\begin{equation}
    F\left(y_{i}, y_{i+1}, t_{i+1}, t_{i}\right) = 0.
    \label{eq:implicit}
\end{equation}
We specialize our results for a specific integration scheme, the backward Euler method defined by the implicit 
function
\begin{equation}
    y_{i+1} - y_{i} - h\left(y_{i+1}, t_{i+1}; p \right) \Delta t_{i+1} = 0
    \label{eq:backward-euler}
\end{equation}
with $\Delta t_{i+1} = t_{i+1} - t_{i}$.

In general, Eq. \ref{eq:implicit} is nonlinear.  Newton's method provides a means to solve the nonlinear system with the iterative scheme
\begin{equation}
    \prescript{}{k+1}{y}_{i+1} = \prescript{}{k}{y}_{i+1} - \prescript{}{k}{J}_{i+1}^{-1}\prescript{}{k}{F}_{i+1}
\end{equation}
with
\begin{equation}
    \prescript{}{k}{J}_{i+1}^{-1} = \frac{d \prescript{}{k}{F}_{i+1}}{d \prescript{}{k}{y}_{i+1}}
    \label{eq:discrete-residual}
\end{equation}
the discrete Jacobian and 
\begin{equation}
    \prescript{}{k}{F}_{i+1} = F\left(y_{i}, \prescript{}{k}{y}_{i+1}, t_{i+1}, t_{i}\right).
    \label{eq:discrete-jacobian}
\end{equation}
Here we introduce yet another subscript, now before a quantity, which indicates the iteration count in Newton's method.  In addition to this  recursion formula, Newton's method requires a starting guess at the converged solution, often taken as the converged value at the previous time step
\begin{equation}
    \prescript{}{0}{y}_{i+1} = y_{i}.
\end{equation}
Again we will omit these subscripts except where necessary.

For the backward Euler method the update formula is
\begin{equation}
    \prescript{}{k+1}{y}_{i+1} = \prescript{}{k}{y}_{i+1} - \prescript{}{k}{J}_{i+1}^{-1}\left[\prescript{}{k}y_{i+1} - y_{i} - h\left(\prescript{}{k}y_{i+1}, t_{i+1}; p \right) \Delta t_{i+1}\right]
\end{equation}
and the discrete Jacobian is 
\begin{equation}
    J_{i+1} = I - j_{i+1} \Delta t_{i+1}
\end{equation}
with $j$ the ODE Jacobian
\begin{equation}
    j_{i+1} = \frac{d h\left(y_{i+1}, t_{i+1} ; p \right)}{d y_{i+1}}.
\end{equation}

\paragraph{Explicit methods}

For comparison, we also reimplement the explicit, forward Euler approach used in \cite{chen2018neural}.  The update equation for this time integration scheme is:
\begin{equation}
    y_{i+1} - y_{i} - h\left(y_{i}, t_{i}; p \right) \Delta t_{i+1} = 0.
    \label{eq:forward-euler}
\end{equation}

This equation is \emph{not} implicit in the next solution $y_{i+1}$ and so the implementation does not have to solve a nonlinear equation to determine the next state.

\subsubsection{Calculating the ODE Jacobian}

The backward Euler method requires both the ODE function $h$ and its derivative $j$.  This work considers three options for calculating $j$:
\begin{enumerate}
    \item an analytic definition
    \item backward mode automatic differentiation
    \item forward mode automatic differentiation.
\end{enumerate}
The analytic derivative of the problem with respect to state $y$ is often much easier to derive than the partial derivative with respect to the problem parameter $p$.  As shown below, this analytic derivative is tractable even for neural ODEs.  However, AD options provide a means to applying implicit integration methods that do not require any additional calculus from the user, at the expense of additional computational work and memory use.

\subsubsection{Gradient calculation}

\citet{chen2018neural} compares the efficiency of calculating the gradient of the integrated response (Eq. \ref{eq:gradient}) using backward mode AD and the adjoint method for explicit, batched time integration.  Simple modifications extend both approaches to implicit time integration schemes.  Appendix \ref{subsec:adjoint-simple} derives the adjoint approach for backward Euler time integration.

\subsection{Chunked time integration}

\subsubsection{Time integration}

The main problem this paper addresses is GPU bandwidth.  If the product of the batch size ($n_{batch}$) and the problem size ($n_{size}$) is comparatively small then the time integration (and adjoint gradient calculation algorithms) will starve the GPU for work.  There is a third potential axis to batch vectorize: time.  As noted above, even for the case where the application only requires final value of the solution trajectory, numerical time integration must often reach this final solution incrementally, integrating through multiple time steps.  Vectorizing the integration and adjoint algorithms along this time axis provides significantly more parallel work to the GPU.

Naively, we might expect parallelizing time integration to be impossible: after all in Eq. \ref{eq:backward-euler} and \ref{eq:forward-euler} each time step depends on the solution of the ODE at the previous step.  However, taking inspiration from space time finite element methods \cite{hughes1988space}, we can indeed vectorize the algorithm in time.  The following presents a derivation of this chunked time integration for the implicit, backward Euler time integration scheme.  Vectorizing the forward Euler scheme is similar and, in fact, more straightforward as it does not require solving an implicit system of equations.  Our implementation includes both time integration schemes, but we describe only the backward Euler formulation in detail here.

Consider a ``chunk'' of $n_{chunk}$ time steps integrating the system of ODEs from time $t_i$ to time $t_{i+n_{chunk}}$.  Recast the solution at each of these time steps as
\begin{equation}
    y_{i+j} = y_{i} + \Delta y_{j}.
\end{equation}
That is, the value at a given time is equal to the value at the start of the chunk of time steps plus some increment.  With this rearrangement, we can write the implicit time integration equations (Eq. \ref{eq:implicit}) for the entire chunk of $n_{chunk}$ steps as
\begin{equation}
    \begin{bmatrix}
    F\left(y_{i}, y_{i} + \Delta y_{1}, t_{i+1}, t_{i}\right)\\ 
    F\left(y_{i}+ \Delta y_{1}, y_{i} + \Delta y_{2}, t_{i+2}, t_{i+1}\right)\\ 
    \vdots\\ 
    F\left(y_{i} + \Delta y_{j}, y_{i} + \Delta y_{j-1}, t_{i+j}, t_{i+j-1}\right)
    \end{bmatrix} = 0.
\end{equation}
Or, specifically for the backward Euler method
\begin{equation}
    \begin{bmatrix}
    \Delta y_{1} - h\left(y_{i} + \Delta y_{1}, t_{i+1}; p \right) \Delta t_{i+1}\\ 
    \Delta y_{2} - \Delta y_{1} - h\left(y_{i} + \Delta y_{2}, t_{i+2}; p \right) \Delta t_{i+2}\\ 
    \vdots\\ 
    \Delta y_{j} - \Delta y_{j-1} - h\left(y_{i} + \Delta y_{j}, t_{i+j}; p \right) \Delta t_{i+j}
    \end{bmatrix} = 0.
    \label{eq:batched-backward-euler}
\end{equation}

In the simple, unchunked time integration scheme (Eq. \ref{eq:backward-euler}) the nonlinear system has batch size $n_{batch}$ and dimension $n_{size}$.  The nonlinear system for this chunked time integration scheme (Eq. \ref{eq:batched-backward-euler}) has batch size $n_{batch}$ and dimension $n_{size} \times n_{chunk}$.  Clearly $n_{chunk} \le n_{time}$ but typically $n_{time}$ is much larger than $n_{size}$ (and often larger than $n_{batch}$) meaning the approach can fully consume the available GPU bandwidth by tuning $n_{chunk}$.

Equation \ref{eq:batched-backward-euler} is still nonlinear.  Applying Newton's method provides the update algorithm
\begin{equation}
    \begin{split}
    \begin{bmatrix}
    \prescript{}{k+1} \Delta y_{1}\\ 
    \prescript{}{k+1} \Delta y_{2}\\ 
    \vdots\\ 
    \prescript{}{k+1} \Delta y_{j}
    \end{bmatrix} = & \begin{bmatrix}
    \prescript{}{k} \Delta y_{1}\\ 
    \prescript{}{k} \Delta y_{2}\\ 
    \vdots\\ 
    \prescript{}{k} \Delta y_{j}
    \end{bmatrix} \\  &- \prescript{}{k}{J}_{i+1}^{-1} \begin{bmatrix}
    \prescript{}{k} \Delta y_{1} - h\left(y_{i} + \prescript{}{k} \Delta y_{1}, t_{i+1}; p \right) \Delta t_{i+1}\\ 
    \prescript{}{k} \Delta y_{2} - \prescript{}{k} \Delta y_{1} - h\left(y_{i} + \prescript{}{k} \Delta y_{2}, t_{i+2}; p \right) \Delta t_{i+2}\\ 
    \vdots\\ 
    \prescript{}{k} \Delta y_{j} - \prescript{}{k} \Delta y_{j-1} - h\left(y_{i} + \prescript{}{k} \Delta y_{j}, t_{i+j}; p \right) \Delta t_{i+j}
    \end{bmatrix}
    \end{split}
    \label{eq:batched-update}
\end{equation}
where the chunked discrete Jacobian has a bidiagonal form
\begin{equation}
    \prescript{}{k}{J}_{i+1} = \begin{bmatrix}
I - \prescript{}{k}j_{i+1} & 0  & \cdots & 0 & 0 \\ 
-I & I - \prescript{}{k}j_{i+2}  & 0  & \cdots & 0 \\ 
0 & -I & I - \prescript{}{k}j_{i+3}  & \cdots  & 0 \\ 
\vdots  & \vdots  & \ddots & \ddots  & \vdots \\
0 & 0 & 0  & -I &  I - \prescript{}{k}j_{i+j}
\end{bmatrix}
\end{equation}
where $\prescript{}{k}j_{i+j}$ is the ODE Jacobian evaluated at the current values of the state $y_i + \prescript{}{k}\Delta y_{j}$.

The key linear algebra kernel for the chunked time integration is solving batched, blocked, bidiagonal systems of this type to update the current Newton-Raphson iterate for the incremental solution.  Section \ref{sec:bidiagonal} discusses solving these systems efficiently on GPUs.

Again, the algorithm also requires a guess for the incremental solution vector $\Delta y$ to start the iterations.  We consider the simplest possible guess
\begin{equation}
    \prescript{}{0}{\Delta y_{i+j}} = 0
\end{equation}
which assumes the values of the ODE do not change from the last time step in the previous chunk.

Algorithm \ref{alg:backward-euler} summarizes the backward Euler chunked time integration scheme.  The examples here use $tol_a = 10^{-8}$ for the absolute tolerance and $tol_r = 10^{-6}$ for the relative tolerance.

\begin{algorithm}
\caption{Chunked backward Euler time integration}\label{alg:backward-euler}
\begin{algorithmic}
\State $i \gets 0$
\While{$i < n_{time}$}
\State $\prescript{}{0}{\Delta y}_{i+b} \gets 0$
\State $\prescript{}{0} r_{i+b} \gets $ residual from Eq. \ref{eq:batched-backward-euler}
\State $\prescript{}{0} n_{i+b} = \left\Vert r_{i+b} \right\Vert $
\State $k \gets 0$
\While{$\prescript{}{k} n_{i+b} > tol_a \forall b$ and $\frac{\prescript{}{k} n_{i+b}}{\prescript{}{0} n_{i+b}} > tol_r \forall b$}
\State $\prescript{}{k+1}{\Delta y_{i+b}} \gets$ update from Eq. \ref{eq:batched-update}
\State $\prescript{}{k+1} r_{i+b} \gets $ residual from Eq. \ref{eq:batched-backward-euler}
\State $k \gets k + 1$
\EndWhile
\State $i \gets i + n_{chunk}$
\EndWhile
\end{algorithmic}
\end{algorithm}

The implementation of the forward Euler integration is much simpler, without the need for a Newton-Raphson loop, and so we omit it here.

\subsubsection{Calculating the ODE Jacobian}

The only difference in calculating the ODE Jacobian $j$ for the chunked time algorithm is that an efficient implementation will batch the calculation over the batch dimension of size $n_{batch}$ and the number of time steps in the chunk $n_{chunk}$.  With this caveat, we consider the same three options: analytical, forward mode AD, and backward mode AD.

\subsubsection{Gradient calculation}

Chunked time integration requires no real changes to the backward mode AD approach for calculating the problem gradient.  The adjoint method extends to chunked/vectorized time integration in much the same manner as forward time integration.  Appendix \ref{subsec:adjoint-chunk} contains a derivation.  Time chunking vectorizes the backward, adjoint march though time in the same way as forward time integration so both the forward and backward pass can be vectorized in time.

\section{Blocked, bidiagonal matrix solvers \label{sec:bidiagonal}}

The key linear algebra kernel for the implicit batched time integration scheme is solving a batched system of blocked, bidiagonal, linear equations of the type
\begin{equation}
\begin{bmatrix}
A_1 &  &  &  & \\ 
B_1 & A_2 &  &  & \\ 
 & B_2 & A_3 &  & \\ 
 &  & \ddots & \ddots & \\ 
 &  &  & B_{n-1} & A_n
\end{bmatrix} \begin{bmatrix}
x_1\\ 
x_2\\ 
x_3\\ 
\vdots\\ 
x_n
\end{bmatrix} = \begin{bmatrix}
y_1\\ 
y_2\\ 
y_3\\ 
\vdots\\ 
y_n
\end{bmatrix}
\end{equation}
We consider two algorithms for solving these types of systems.  In the following discussion we consider matrices with batch size $n_{batch}$, block size $n_{size}$, and square matrix sizes of $n_{chunk}$ blocks in each dimension.  The unrolled array is then a batched matrix of size $n_{batch} \times \left( n_{size} n_{nchunk} \right) \times \left( n_{size} n_{nchunk} \right)$.

\subsection{Thomas's algorithm}

Thomas's algorithm sweeps down the diagonal of the matrix, solving for each (batched, blocked) unknown value $x_i$ one-by-one.  Algorithm \ref{alg:thomas} summarizes the process.  All the matrix operations in this algorithm are batched over an implicit batch dimension, with size $n_{batch}$.

\begin{algorithm}
\caption{Thomas's algorithm for solving batched, blocked matrix systems.}\label{alg:thomas}
\begin{algorithmic}
\State $x_1 = A_{1}^{-1} y_1$
\State $i \gets 1$
\While{$i < n_{chunk}$}
\State $x_{i+1} \gets A_{i+1}^{-1} \left(y_{i+1} - B_{i} x_{i} \right)$
\State $i \gets i + 1$
\EndWhile
\end{algorithmic}
\end{algorithm}

Without considering parallelization/vectorization over the batched matrix operations this algorithm is serial with time complexity 
\begin{equation}
    O(n_{chunk} n_{batch} n_{nblock}^2 + n_{chunk} n_{batch} n_{size}^3) \rightarrow O(n_{chunk} n_{batch} n_{size}^3)
    \label{eq:thomas}
\end{equation} 
with the first cost representing the back solve loop and the second the cost of factorizing the diagonal blocks.  It will be useful to consider the cost of the back solve separately from the cost to factorize the diagonal blocks when comparing to our second algorithm, below, even though asymptotically the factorization cost governs.  Note this is considerably cheaper than solving the full dense matrix, which would be $O(n_{chunk}^3 n_{batch} n_{size}^3)$.

\subsection{Parallel cyclic reduction (PCR) \label{sec:pcr}}

Parallel cyclic reduction is a divide-and-conquer algorithm that recursively splits a (batched, blocked) bidiagonal matrix into two independent bidiagonal linear systems.  The following equations give the basic update, splitting a single bidiagonal system 
\begin{equation}
\begin{bmatrix}
A_1 & 0 & 0 & 0\\ 
B_1 & A_2 & 0 & 0\\ 
0 & B_2 & A_3 & 0\\ 
0 & 0 & B_3 & A_4
\end{bmatrix} \begin{bmatrix}
x_1\\ 
x_2\\ 
x_3\\ 
x_4
\end{bmatrix} = \begin{bmatrix}
y_1\\ 
y_2\\ 
y_3\\ 
y_4
\end{bmatrix}
\end{equation}
into two independent bidiagonal systems
\begin{equation}
\begin{bmatrix}
{\color{red} A_1} & 0 & {\color{red} 0} & 0\\ 
0 & {\color{green} A_2} & 0 & {\color{green} 0}\\ 
{\color{red} -B_2 A_2^{-1} B_1} & 0 & {\color{red} A_3} & 0\\ 
0 & {\color{green} -B_3 A_3^{-1} B_2} & 0 & {\color{green} A_4}
\end{bmatrix} \begin{bmatrix}
{\color{red} x_1}\\ 
{\color{green} x_2}\\ 
{\color{red} x_3}\\ 
{\color{green} x_4}
\end{bmatrix} = \begin{bmatrix}
{\color{red} y_1}\\ 
{\color{green} y_2 - B_1 A_1^{-1} y_1}\\ 
{\color{red} y_3 - B_2 A_2^{-1} y_2}\\ 
{\color{green} y_4 - B_3 A_3^{-1} y_3}
\end{bmatrix}
\end{equation}
where the colors indicate the two independent systems.

PCR then recursively applies this splitting formula to further subdivide the matrices until it reaches a block diagonal form.  Algorithm \ref{alg:pcr} describes the complete process.  The key point in PCR is that the new matrices produced by each application of the splitting formula can be factored independently, in parallel.

\begin{algorithm}
\caption{Parallel cyclic reduction algorithm for solving batched, blocked matrix systems.}\label{alg:pcr}
\begin{algorithmic}
\State $i \gets 1$
\While{$i<\log_2{n_{chunk}}$}
\State Form $2^i$ submatrices by taking every $2^i$th row
\For{submatrix $j=1$ to $j=2^i$} \Comment{This loop can be parallelized}
\For{index $k\in$ submatrix $j$}
\State $B_k \gets -B_k D_k^{-1}B_{k-1}$
\State $y_k \gets y_k - B_k D_{k-1}^{-1} y_{k-1}$
\EndFor
\EndFor
\State $i \gets i + 1$
\EndWhile \\
\Return $x \gets y$
\end{algorithmic}
\end{algorithm}

Without accounting for the parallelism in processing the recursive submatrices, assuming we store the factorization of the $A^{-1}$ blocks, and not accounting for any paralleliztion of the matrix operations over the batch dimension, the time complexity of PCR is
\begin{equation}
    O(\log_2(n_{chunk}) n_{chunk} n_{batch} n_{size}^3 + n_{batch} n_{chunk} n_{size}^3)
\end{equation}
with the first cost from the recursive reduction and the second from the need to first factorize the diagonal blocks.  If our execution device can accommodate the full amount of available parallel work ($n_{chunk}$ submatrices on the final iteration), which is a reasonable assumption for a GPU and the size of matrix systems we consider here, then the time complexity becomes 
\begin{equation}
    O(\log_2(n_{chunk}) n_{batch} n_{size}^3 + n_{batch} n_{chunk} n_{size}^3) \rightarrow O(n_{batch} n_{chunk} n_{size}^3)
\end{equation}
again with the reduction cost first and the factorization cost second.

The two approaches therefore have the same asymptotic performance dominated by the cost of factorizing the diagonal blocks, assuming perfect parallelization of PCR.  However, in practice, as explored in the timing studies below, the two algorithms are competitive depending on the size of the system (defined by $n_{chunk}$, $n_{size}$, and $n_{batch}$), as suggested by the differing complexity of the reduction/backsolve part of the analysis.

\subsection{A hybrid approach}

A third option is to start solving the system with PCR, but halt the recursive subdivision before reducing the system to a series of diagonal blocks and instead solve the remaining, unreduced sets of independent equations using Thomas's algorithm.  This approach has an extra parameter controlling the heuristic, the number of iterations after which to switch from PCR to Thomas, $n_{switch}$.

Our implementation includes this hybrid approach, but we found it did not outperform the straight Thomas's algorithm or PCR for any problem and any value of $n_{switch}$ we evaluated.  As such, we do not include the results for the hybrid approach in our discussion below.

\section{Notes on the implementation \label{sec:implementation}}

We implement chunked time integration with both backward and forward Euler approaches in the open source \emph{pyoptmat} package, available at \url{https://github.com/Argonne-National-Laboratory/pyoptmat}.  \emph{pyoptmat} uses pytorch \cite{pytorch} to provide efficient tensor operations on the GPU.  The timing results below were generated with \texttt{v1.3.5}, available as a tag in the repository.  The example problems, detailed below, are available in a separate repository, available at \url{https://github.com/reverendbedford/odeperf}.

The following subsections contain a few notes on the implementation.

\subsection{Chunked time integration}

\emph{pytorch} includes additional modules specializing the general ODE integration routines discussed here to infer material constitutive models from data.  The \texttt{ode} and \texttt{chunktime} modules contain the routines described and used in this work.  The key routine in the \texttt{ode} module is \texttt{ode.odeint\_adjoint} which implements chunked forward and backward Euler time integration for generic systems of ordinary differential equations, encoded as pytorch \texttt{nn.Module} classes.

The only annoyance for the end users is that their ODE models must accept two batch dimensions -- one over the number of time series to consider at once ($n_{batch}$) and a second over the number of time steps to consider at once ($n_{chunk}$).  Writing most models in this way is usually trivial.

This routine provides pytorch gradients using the adjoint method.  A related routine \texttt{ode.odeint} provides gradients via pytorch's backward mode AD implementation, for comparison.

\subsection{Calculating the Jacobian}

Implicit schemes require both the rate of the ODE system (Eq. \ref{eq:discrete-residual}) and also the Jacobian of this system with respect to the state (Eq. \ref{eq:discrete-jacobian}).  As noted above, we compare the cost of implicit time integration calculating this Jacobian three ways: analytically, with backward mode AD, and with forward mode AD.

For backward and forward mode AD we vectorize the Jacobian calculation over the batch and chunked time dimensions with the \texttt{functorch} tools recently incorporated into the main pytorch library.

\subsection{Batched, blocked bidagonal solvers}

The implementation of Thomas's algorithm in pytorch/\emph{pyoptmat} is straightforward.

PCR is less straightforward.  pytorch only provides parallelism via tensor axis vetorization, not, for example, via threaded parallelism.  Most discussion on PCR provides a recursive, threaded implementation which we cannot implement efficiently in pytorch.  The iterative description in Algorithm \ref{alg:pcr} matches, in principle, the actual implementation in \emph{pyoptmat}.  However, even here it may not be clear how to parallelize the submatrix loop noted in the algorithm description.

Our solution is to provide a view into the matrix storage with an extra axis representing the current set of submatrices.  Our implementation initially stores the matrix $A$ as a tensor with dimension $(n_{chunk}$,$n_{batch}$,$n_{size}$,$n_{size})$.  To implement PCR we reshape this to be a tensor of size $(1,n_{chunk},n_{batch},n_{size},n_{size})$.  Subsequent iterations in the parallelizable loop ``shuffle'' the first two indices to provide submatrix views; first reshaping the data to a shape of $(2$,\allowbreak$n_{chunk}/2$,\allowbreak$n_{batch}$,\allowbreak$n_{size}$,\allowbreak$n_{size})$, then to $(4$,$n_{chunk}/4$,$n_{batch}$,$n_{size}$,$n_{size})$, etc.  It is possible to implement these tensor reshape operations without copying the data by manipulating the tensor strides.  Only by doing this could we implement a version of PCR that can efficiently compete with Thomas's algorithm.

A second issue is how to deal with matrices with $n_{chunk}$ not equal to some power of two.  One option is a variant of the hybrid algorithm briefly described above: using PCR for the largest power of two sized submatrix and Thomas for the remainder.  Instead, we implement an iterative PCR approach which factorizes the matrix size into a sum of powers of two (i.e. $n_{chunk} = 7 = 4 + 2 + 1$) and apply PCR iteratively to each power of two-sized submatrix in turn.

\section{Examples \label{sec:example-sec}}

\subsection{Example problems \label{subsec:examples}}

We completed a timing study on several representative systems of differential equations to assess the performance of the implicit, vertorized time integration algorithm.  The sample ODEs span a range of applications, from simple to more complex physical systems to neural ODEs.  These ODEs sample a reasonable range of problems, capturing the performance of the implicit, chunked time integration approach developed here for a wider variety of applications.  We used sample systems where we can arbitrarily tune the size of the system ($n_{size}$) and the total number of time steps ($n_{time}$) to explore the effectiveness of the integration algorithm for various size problems.  For the physical problems, the size of the ODE system often must conform to some multiple of a basic unit size.  For these problems we define $n_{size}$ in terms of this $n_{unit}$ dimension.

\subsubsection{Sample ODEs}

\paragraph{Mass-damper-spring system}

\begin{figure}
    \centering
    \includegraphics{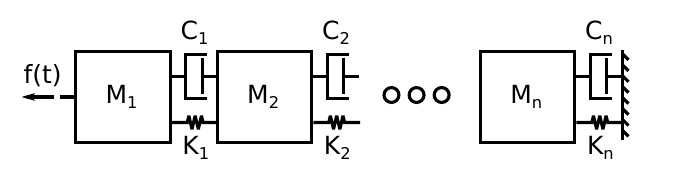}
    \caption{Schematic of the mass-damper-spring system.}
    \label{fig:mds-sketch}
\end{figure}

This system describes a coupled mass-damper-spring system consisting of $n_{unit}$ individual mass-damper-spring units hooked in series, pinned to zero displacement on one end of the chain, and with a periodic applied force $f(t)$ on the other end of the chain (see Figure \ref{fig:mds-sketch}).  To define the system in terms of the balance of linear momentum we must solve for both the displacement $d_i$ and velocity $v_i$ of each element, meaning $n_{size} = 2 n_{unit}$.  We impose zero initial conditions for both the displacements and the velocities.

This problem can be defined by the system of ODEs defining the rate of change of velocity and displacement for each element in the chain
\begin{align} 
    \dot{d}_i &=  v_i \\ 
    \dot{v}_i &=  \frac{K_i}{M_i}\left(d_{i} - d_{i-1} \right) - \frac{K_{i+1}}{M_{i+1}}\left(d_{i+1} - d_{i} \right) + \frac{C_i}{M_i}\left(v_{i} - v_{i-1} \right) \\ \nonumber & - \frac{C_{i+1}}{M_{i+1}}\left(v_{i+1} - v_{i} \right) + f_i
\end{align}
with the convention that $u_{n_{unit}+1} = v_{n_{unit}+1} = 0$,
\begin{equation}
f_i = \left\{\begin{matrix}
f(t) & i = 1 \\ 
0 & i > 1
\end{matrix}\right.
\end{equation}
and
\begin{equation}
    f(t) = f_{a} \sin\left(\frac{2\pi}{T}\right).
\end{equation}
The convention here is element $i=1$ has the force applied and element $i=n_{unit}+1$ is fixed.

Table \ref{tab:mass-params} provides the particular parameters used in the numerical study.

\begin{table}
\centering
\begin{tabular}{ll}
\hline
Parameter & Values                              \\ \hline
$K_i$     & \texttt{linspace($10^{-2}$,$10^{0}$,$n_{unit}$)} \\
$C_i$     & \texttt{linspace($10^{-6}$,$10^{-4}$,$n_{unit}$)} \\
$M_i$     & \texttt{linspace($10^{-7}$,$10^{-5}$,$n_{unit}$)} \\
$f_a$     & 1.0                                 \\
$T$       & \texttt{linspace($10^{-2}$,$10^{0}$,$n_{batch}$)} \\
$t_{max}$ & 1.0                                   \\ \hline
\end{tabular}
\caption{Model parameters used in the timing study for the mass-damper-spring case.}
\label{tab:mass-params}
\end{table}

\begin{figure}[!htb]
    \centering
    \begin{subfigure}{0.48\textwidth}
        \centering
        \includegraphics[width=\textwidth]{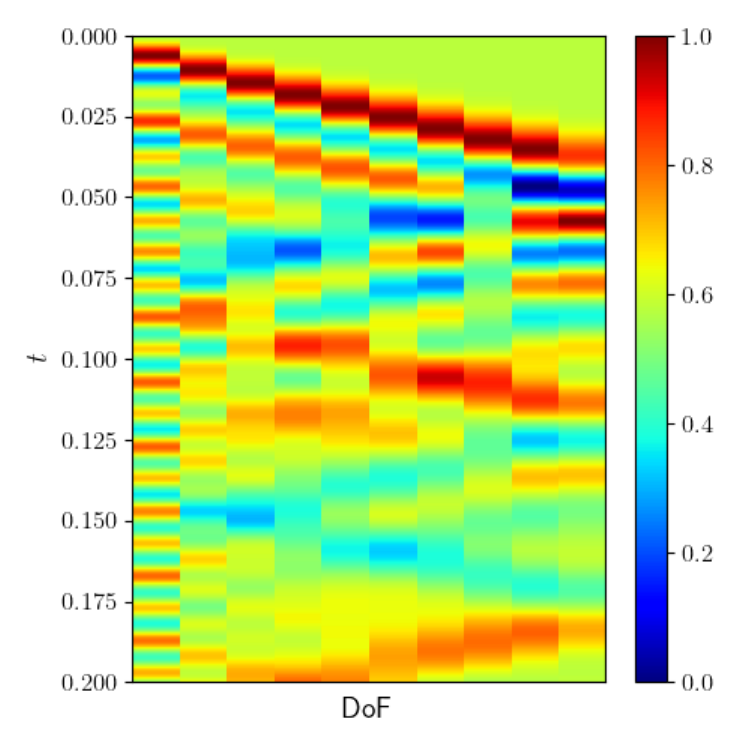}
        \caption{}
    \end{subfigure}
    \begin{subfigure}{0.48\textwidth}
        \centering
        \includegraphics[width=\textwidth]{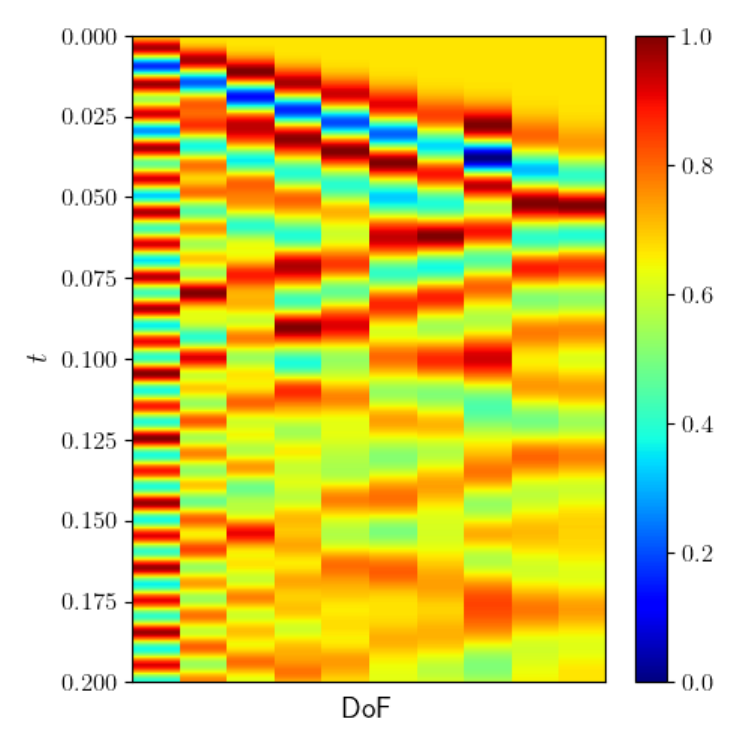}
        \caption{}
    \end{subfigure}
    \caption{Sample normalized trajectory for the mass-damper-spring problem with $n_{size}=10$ and $t_{max}=0.2$. All the other controllable parameters are using their respective lower bounds. (a) Normalized displacements. (b) Normalized velocities.}
    \label{fig:sample-mds}
\end{figure}

Figure \ref{fig:sample-mds} shows a sample trajectory.


\paragraph{Coupled neuron model}

This model, detailed by \citet{schwemmer2012theory}, describes the response of a collection of weakly-coupled neurons to an oscillating current.  The ODE model describes each neuron in terms of four key variables: the voltage, $V_i$, and three non-dimensional probability values, $m$, $h$, and $n$, related to the potassium channel activation, sodium channel activation, and sodium channel deactivation, respectively.  Four ODEs represent the response of each neuron, so the size of this model is $n_{size} = 4 n_{unit}$.

The equations representing the response of a single neuron are:
\begin{align}
  \phantom{\dot{V}_i}
  &\begin{aligned}
    \mathllap{\dot{V}_i} &= \frac{1}{C_i} \biggl[-g_{Na,i} m_i^3 h_i \left(V_i - E_{Na,i} \right) \\
      &\qquad -g_{K,i} n_i^4 \left(V_i - E_{K,i} \right) -g_{L,i} \left(V_i - E_{L,i}\right) \\
      &\qquad + I_i(t) + g_{C,i} \sum_{j=1}^{n_{unit}} \left(V_i - V_j \right) \biggr]
  \end{aligned}\\
  &\begin{aligned}
    \mathllap{\dot{m}_i} &=  \frac{m_{\infty,i} - m_i}{\tau_{m,i}}
  \end{aligned}\\
  &\begin{aligned}
    \mathllap{\dot{h}_i} &= \frac{h_{\infty,i} - h_i}{\tau_{h,i}}
  \end{aligned}\\
  &\begin{aligned}
    \mathllap{\dot{n}_i} &= \frac{n_{\infty,i} - n_i}{\tau_{n,i}}
  \end{aligned}
\end{align}
with
\begin{equation}
    I_i(t) = I_{a} \sin\left(\frac{2 \pi}{T} t \right).
\end{equation}

Table \ref{tab:neuron-params} provides the parameters used in the numerical study.

\begin{table}
\centering
\begin{tabular}{ll}
\hline
Parameter & Values                              \\ \hline
$C_i$     & \texttt{linspace($10^{-1}$,$10^{0}$,$n_{unit}$)} \\
$g_{Na,i}$     & \texttt{linspace($10^{-1}$,$10^{0}$,$n_{unit}$)} \\
$E_{Na,i}$     & \texttt{linspace($10^{-1}$,$10^{0}$,$n_{unit}$)} \\
$g_{K,i}$     & \texttt{linspace($10^{-1}$,$10^{0}$,$n_{unit}$)}\\
$E_{K,i}$       & \texttt{linspace($10^{-1}$,$10^{0}$,$n_{unit}$)} \\
$g_{L,i}$ & \texttt{linspace($10^{-1}$,$10^{0}$,$n_{unit}$)} \\
$E_{L,i}$ & \texttt{linspace($10^{-1}$,$10^{0}$,$n_{unit}$)}\\
$m_{\infty,i}$ & \texttt{linspace($10^{-1}$,$10^{0}$,$n_{unit}$)} \\
$\tau_{m,i}$ & \texttt{linspace($0.5$,$5$,$n_{unit}$)} \\ 
$h_{\infty,i}$ & \texttt{linspace($10^{-1}$,$10^{0}$,$n_{unit}$)} \\
$\tau_h$ & \texttt{linspace($1.5$,$15$,$n_{unit}$)} \\ 
$n_{\infty,i}$ & \texttt{linspace($10^{-1}$,$10^{0}$,$n_{unit}$)} \\
$\tau_{n,i}$ & \texttt{linspace($10^{0}$,$10^{1}$,$n_{unit}$)} \\ 
$g_{C,i}$ & \texttt{linspace($10^{-3}$,$10^{-2}$,$n_{unit}$)} \\
$I_{a,i}$ & \texttt{linspace($10^{-1}$,$10^{0}$,$n_{batch}$)} \\ 
$T_{i}$ & \texttt{linspace($0.5$,$2$,$n_{unit}$)} \\
$t_{max}$ & 10.0                                \\ \hline
\end{tabular}
\caption{Model parameters used in the timing study for the coupled neuron case.}
\label{tab:neuron-params}
\end{table}

\begin{figure}[!htb]
    \centering
    \begin{subfigure}{0.48\textwidth}
        \centering
        \includegraphics[width=\textwidth]{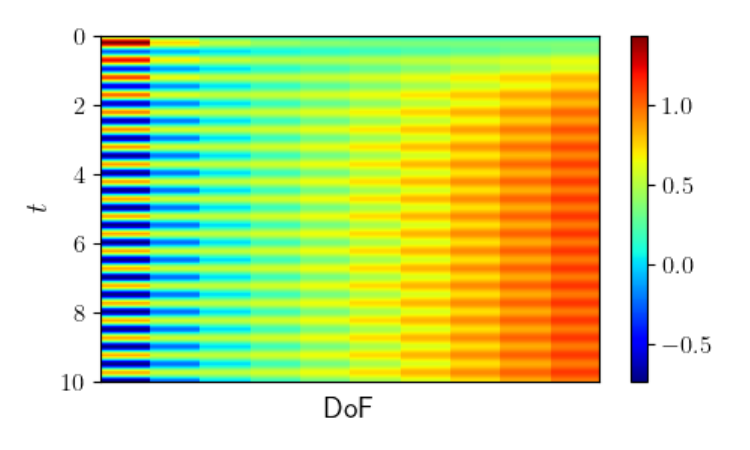}
        \caption{}
    \end{subfigure}
    \begin{subfigure}{0.48\textwidth}
        \centering
        \includegraphics[width=\textwidth]{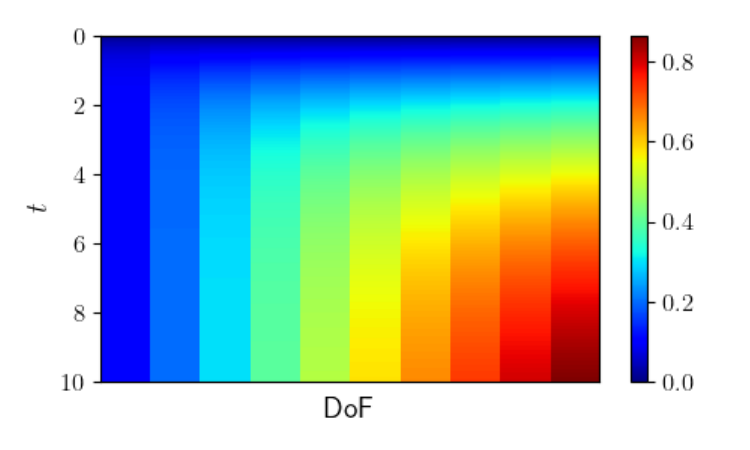}
        \caption{}
    \end{subfigure}
    \begin{subfigure}{0.48\textwidth}
        \centering
        \includegraphics[width=\textwidth]{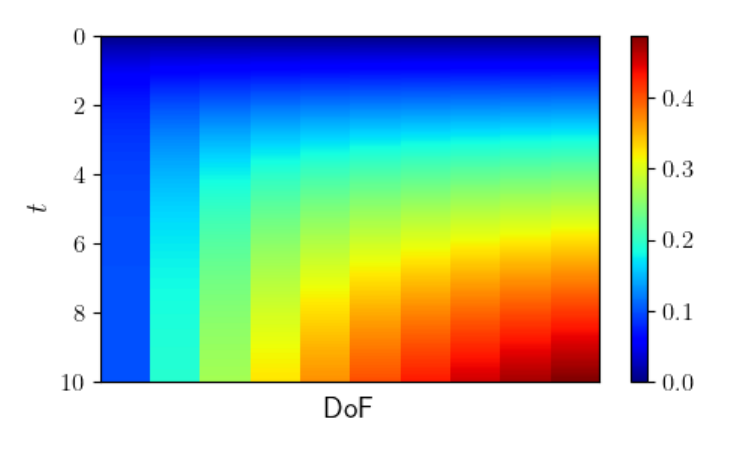}
        \caption{}
    \end{subfigure}
    \begin{subfigure}{0.48\textwidth}
        \centering
        \includegraphics[width=\textwidth]{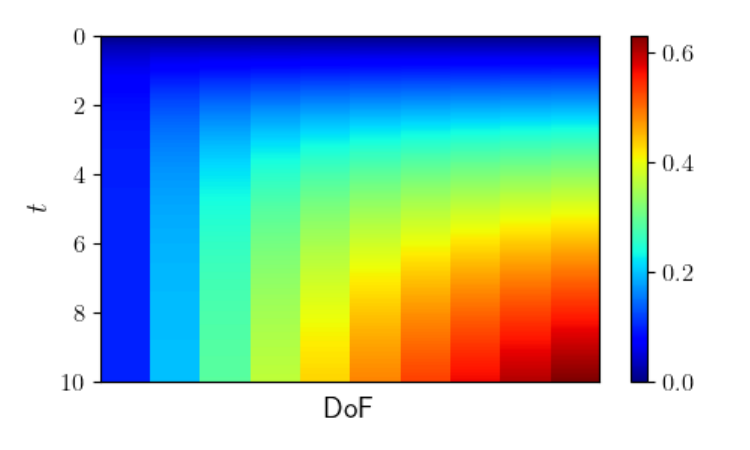}
        \caption{}
    \end{subfigure}
    \caption{Sample trajectory for the neuron problem with $n_{size}=10$, $t_{max}=10$, and $T=0.5$. All the other controllable parameters are using their respective upper bounds. (a) Voltage. (b-d) Non-dimensional probabilities $m$, $h$ and $n$.}
    \label{fig:sample-neuron}
\end{figure}

Figure \ref{fig:sample-neuron} shows a sample trajectory.


\paragraph{Chaboche viscoplastic constitutive model}

This model represents the response of a metallic material under high temperature, cyclic loading \cite{chaboche1989constitutive}.  The model includes one variable representing the evolution of stress in the material, one describing the material strengthening under monotonic load, and an arbitrary number describing the material strengthening under cyclic load.  As such $n_{size} = 2 + n_{unit}$.

The equations describing the system of ODEs are:
\begin{equation}
    \dot{y} = \begin{bmatrix}
\dot{\sigma}\\ 
\dot{K}\\ 
\dot{X}_1\\ 
\vdots \\ 
\dot{X}_{n_{unit}}
\end{bmatrix}
\end{equation}
with
\begin{align} 
    \dot{\sigma} &=  E \left(\dot{\varepsilon}(t) - \dot{\varepsilon}_p \right) \\ 
    \dot{K} &=  \tau \left(K_\infty - K \right) \\
    \dot{X}_i &= \frac{2}{3}C_i \dot{\varepsilon}_p - \gamma_i X_i \left| \dot{\varepsilon}_p \right|
\end{align}
with $\langle\rangle$ the Macaulay brackets, 
\begin{equation}
    \dot{\varepsilon}_p = \left\langle \frac{\left| \sigma - \sum_{i=1}^{n_{unit}} X_i \right| - K - \sigma_0}{\eta} \right\rangle^n \operatorname{sign}\left(\sigma - \sum_{i=1}^{n_{unit}} X_i \right),
\end{equation}
and
\begin{equation}
    \dot{\varepsilon}(t) = \dot{\varepsilon}_a \sin\left(\frac{2 \pi}{T} t \right).
\end{equation}

Table \ref{tab:chaboche-params} describes the parameters used in the scaling study.

\begin{table}
\centering
\begin{tabular}{ll}
\hline
Parameter & Values                              \\ \hline
$E$     & 10.0 \\
$n$     & 5.0 \\
$\eta$     & 2.0 \\
$\sigma_0$     & 1.0                                 \\
$K_\infty$ & 10.0 \\
$\tau$ & 1.0 \\
$C_i$ & \texttt{linspace($0.1$,$1.0$,$n_{unit}$)} \\
$\gamma_i$ & \texttt{linspace($0.1$,$0.5$,$n_{unit}$)} \\
$\dot{\varepsilon}_a$ & \texttt{linspace($10^{-1}$,$10^{0}$,$n_{batch}$)} \\
$T$       & 1.0 \\
$t_{max}$ & 10.0                                   \\ \hline
\end{tabular}
\caption{Model parameters used in the timing study for the Chaboche model case.}
\label{tab:chaboche-params}
\end{table}

\begin{figure}[!htb]
    \centering
    \begin{subfigure}{0.48\textwidth}
        \centering
        \includegraphics[width=\textwidth]{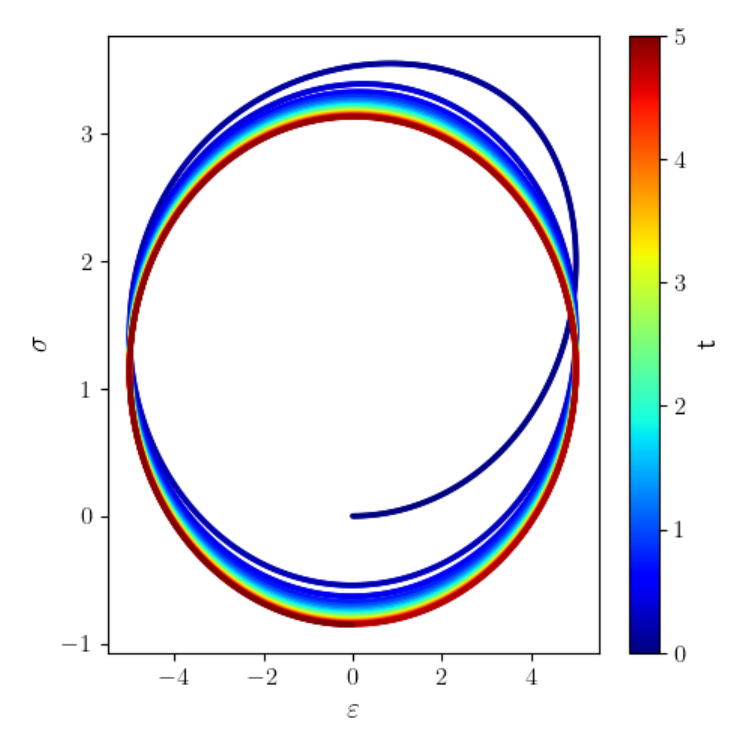}
        \caption{}
    \end{subfigure}
    \begin{subfigure}{0.48\textwidth}
        \centering
        \includegraphics[width=\textwidth]{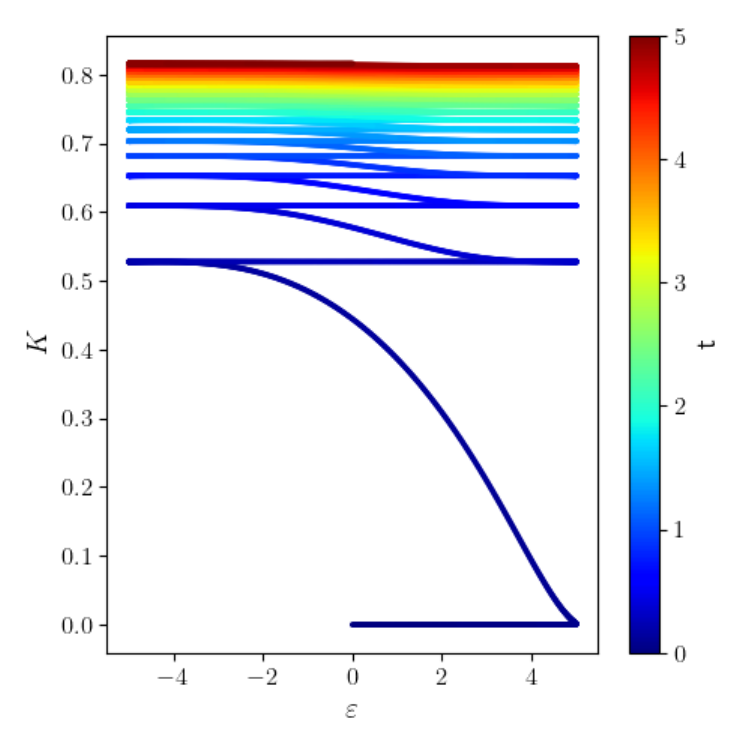}
        \caption{}
    \end{subfigure}
    \caption{Sample trajectory for the Chaboche problem with 10 backstresses, $t_{max}=5$, and $T=0.25$. All the other controllable parameters are using their respective upper bounds. (a) The stress-strain hysterisis loop (b) The hardening as a function of strain. Both plots are colored by the marching time $t$.}
    \label{fig:sample-chaboche}
\end{figure}

Figure \ref{fig:sample-chaboche} shows a sample trajectory.


\paragraph{Neural ODE}

The final example problem is a neural ODE driven by a periodic forcing function.  The equations defining this model are
\begin{equation}
    \dot{y} = \tanh \circ L_3 \circ \tanh \circ L_2 \circ \tanh \circ L_1 \circ \begin{bmatrix}
y\\ 
f(t)
\end{bmatrix}
\end{equation}
where $y$ has size $n_{unit}$, $L_1$, $L_2$, and $L_3$ are linear layers with $(\mathrm{input},\mathrm{output})$ sizes $(n_{unit}+1,n_{unit}+1)$, $(n_{unit}+1,n_{unit}+1)$, and $(n_{unit}+1,n_{unit})$, respectively, and the driving force is
\begin{equation}
    f(t) = f_a \sin \left( \frac{2 \pi}{T} t \right).
\end{equation}
The linear layer weights are randomly drawn from a uniform distribution in the range $[-\sqrt{\frac{1}{n_{unit} + 1}}, \sqrt{\frac{1}{n_{unit} + 1}}]$, the biases are randomly drawn from a uniform distribution in the same range, and Table \ref{tab:nn-params} provides the forcing function parameters used in the scaling study.  The random weights and biases are fixed for all trials in the numerical studies below.

\begin{table}
\centering
\begin{tabular}{ll}
\hline
Parameter & Values                              \\ \hline
$f_a$ & 1.0 \\
$T$       & \texttt{linspace($10^{-2}$,$10^{0}$,$n_{batch}$)} \\
$t_{max}$ & 1.0                                   \\ \hline
\end{tabular}
\caption{Model parameters used in the timing study for the neural ODE case.}
\label{tab:nn-params}
\end{table}

This problem size has $n_{size} = n_{unit}$.

\subsubsection{Timing study}

The timing results below fix the ODE model parameters but vary the numerical parameters describing the size of the problem, the number of batches, the time integration chunk size, the method for calculating the Jacobian matrix, the integration method, the method for calculating the sensitivity of the model response with respect to the ODE parameters, and the batched, blocked bidiagonal solver.  Table \ref{tab:numerical-params} summarizes the available numerical parameters.  Sections \ref{subsec:adadj}, \ref{subsec:chunk}, and \ref{subsec:adjac} provide the actual values for each parameter study.

\begin{table}
\centering
\begin{tabular}{lll}
\hline
Parameter          & Description                                     & Values                        \\ \hline
$n_{unit}$         & \makecell[l]{Fundamental size, related \\ to the system size, $n_{size}$}   & $\ge1$                        \\
$n_{batch}$        & Number of parallel batches                      & $\ge1$                        \\
$n_{time}$         & Number of time steps                            & $\ge1$                        \\
$n_{chunk}$        & Number of vectorized time steps                 & $1 \le n_{chunk} \le n_{time}$              \\
Jacobian      & Method to calculate the Jacobian                & \makecell[l]{analytic, forward, \\backward}   \\
Backward      & \makecell[l]{Method to calculate \\ the parameter sensitivities} & adjoint, AD                   \\
Solver & Method to solve linearized system & Thomas, PCR \\
Integration & Numerical time integration method & \makecell[l]{backward Euler,\\ forward Euler} \\ \hline
\end{tabular}
\caption{Numerical parameters varied in the timing studies.}
\label{tab:numerical-params}
\end{table}

The actual timing run for each ODE and set of numerical parameters is a forward integration of the model for time $t \in [0,t_{max}]$ using $n_{time}$ equally-sized time steps to provide a time series of $y_i$ values, calculating an arbitrary loss function over the results
\begin{equation}
    L = \sqrt{\sum_{k=1}^{n_{size}}\sum_{i=1}^{n_{time}}\sum_{j=1}^{n_{batch}} \left(y_{i,k}^{(j)}\right)^2},
\end{equation}
i.e. the Frobenius norm of the tensor containing the results for all times and all batches, and calculating the derivative $\frac{dL}{dp}$ where $p$ is the complete set of model parameters with a backward pass.  During this process each study monitors the wall time for the forward pass, the backward pass, and the total wall time and the maximum GPU memory required during any point in both the forward and backward passes.  We repeat each individual trial three times and average the results.  We ran the timing studies on a single NVIDIA RTX A5000 GPU using double precision floating point numbers.

\subsection{AD versus adjoint sensitivity \label{subsec:adadj}}

This is a small study to reconfirm the results of \cite{chen2018neural} for implicit time integration --- namely, that the adjoint approach significantly outperforms automatic differentiation for calculating parameter sensitivities in ODE problems.  Table \ref{tab:ad-v-adjoint} details the parameters varied in the study, which considered all combinations of the lists in the table.  This study focuses on the efficiency of backward mode AD and the adjoint method for calculating parameter sensitivities for functions of the results of ODE time integration as the number of time steps increases, i.e. as the problem gets deeper.

\begin{table}
\centering
{\tiny
\begin{tabular}{l|llllllll}
           & $n_{unit}$ & $n_{batch}$ & $n_{time}$                     & $n_{chunk}$ & Jacobian & Backward & Solver & Integration \\ \hline
\makecell[l]{Neural \\ODE} & 20         & 100         & \makecell[l]{100,200,300,\\500,1000,\\1500,2000} & 50          & analytic      & AD, adjoint & Thomas   & backward    
\end{tabular}
}
\caption{Numerical parameter study comparing AD versus adjoint sensitivity calculations.}
\label{tab:ad-v-adjoint}
\end{table}

\begin{figure}
\centering
\begin{subfigure}{.5\textwidth}
  \centering
  \includegraphics[width=.95\linewidth]{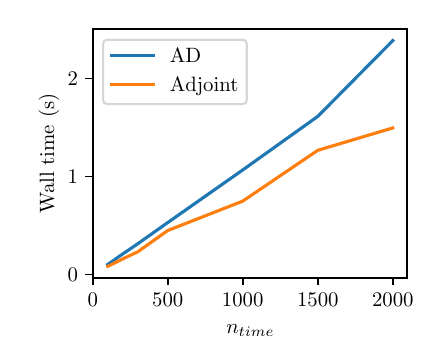}
  \caption{Wall time}
\end{subfigure}%
\begin{subfigure}{.5\textwidth}
  \centering
  \includegraphics[width=.95\linewidth]{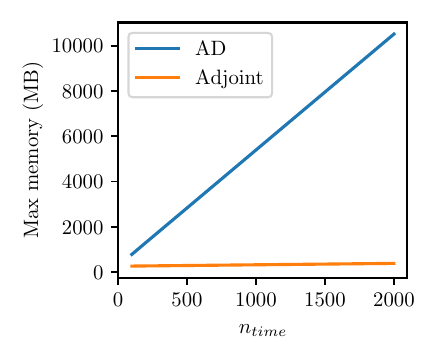}
  \caption{Memory}
\end{subfigure}
\caption{Wall time (a) and memory (b) comparison for the neural ODE case illustrating the effectiveness of the adjoint method in computing parameter sensitivities of functions of the integrated results.}
\label{fig:ad_adjoint}
\end{figure}

Figure \ref{fig:ad_adjoint} plots the wall time and memory efficiency of the two approaches as a function of increasing $n_{time}$.  The adjoint method is slightly faster, but much more memory efficient when compared to AD.  This matches the results reported in \cite{chen2018neural} but now for implicit, backward Euler time integration.  Memory scaling for the adjoint method is linear, like AD, but the slope of the scaling relations is small --- allowing the solver to back propagate through long chains of results from implicit time integration.

This result generalizes to all the ODEs considered here, so in subsequent studies we consider only the adjoint method so we can include longer/larger problems in the performance studies.

\subsection{Efficiency of chunked time integration \label{subsec:chunk}}

Table \ref{tab:big} lists the parameter combinations in this study, focusing on the efficiency of the chunked/vectorized time integration algorithm described here.  This study includes all four types of ODEs outlined in Section \ref{subsec:examples}.

\begin{table}
\centering
{\tiny
\begin{tabular}{l|llllllll}
                   & $n_{unit}$           & $n_{batch}$ & $n_{time}$ & $n_{chunk}$            & Jacobian & Backward & Solver & Integration \\ \hline
\makecell[l]{Mass-\\damper-\\spring} & \makecell[l]{1,2,3,\\5,10}           & \makecell[l]{3,10,\\30,100} & 2000       & \makecell[l]{1,3,10,\\30,100,\\300,1000} & analytic & adjoint  & \makecell[l]{Thomas,\\PCR} & backward    \\
Neuron             & \makecell[l]{1,2,3,\\4,5,6,\\7,8,9,\\10} & \makecell[l]{3,10,\\30,50}  & "          & "                      & "        & "        & "         & "  \\
Chaboche           & \makecell[l]{1,2,3,\\4,5}            & \makecell[l]{3,10,\\30,50}  & "          & "                      & "        & "        & "          & " \\
\makecell[l]{Neural \\ ODE}         & \makecell[l]{1,5,10,\\15,20,25}      & \makecell[l]{3,10,\\30,100} & "          & "                      & "        & "        & "       & "
\end{tabular}
}
\caption{Parameter study for assessing the efficiency of the chunked time integration approach and for comparing the cost of calculating the Jacobian matrix for implicit integration with different methods.}
\label{tab:big}
\end{table}

\begin{figure}
    \centering
    \emph{Mass-damper-spring}
    \begin{subfigure}{1.0\textwidth}
        \centering
        \begin{subfigure}{.375\textwidth}
            \centering
            \includegraphics[width=.95\linewidth]{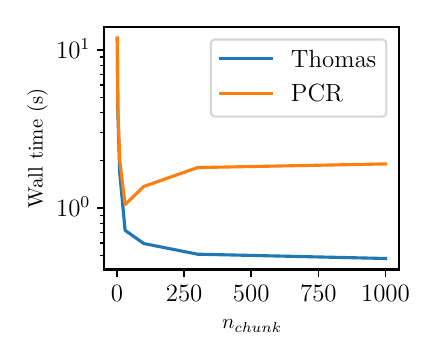}
            \caption{Wall time}
        \end{subfigure}%
        \begin{subfigure}{.375\textwidth}
            \centering
            \includegraphics[width=.95\linewidth]{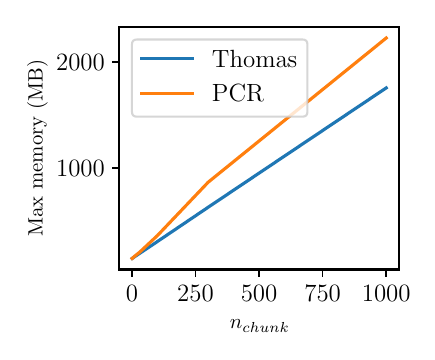}
            \caption{Memory}
        \end{subfigure}
    \end{subfigure}

    \emph{Neuron}
    \begin{subfigure}{1.0\textwidth}
        \centering
        \begin{subfigure}{.375\textwidth}
            \centering
            \includegraphics[width=.95\linewidth]{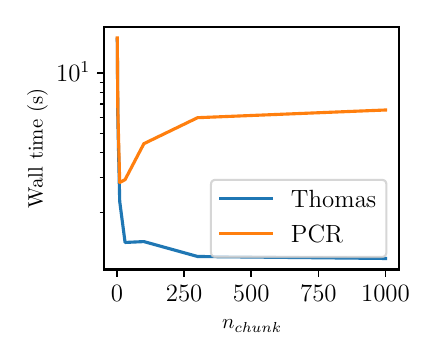}
            \caption{Wall time}
        \end{subfigure}%
        \begin{subfigure}{.375\textwidth}
            \centering
            \includegraphics[width=.95\linewidth]{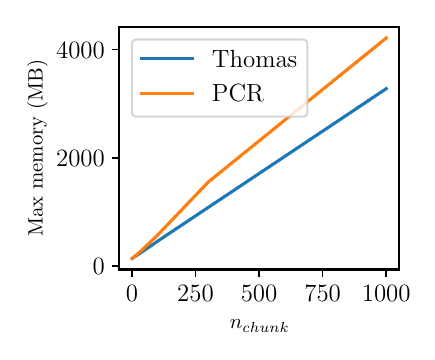}
            \caption{Memory}
        \end{subfigure}
    \end{subfigure}

    \emph{Chaboche}
    \begin{subfigure}{1.0\textwidth}
        \centering
        \begin{subfigure}{.375\textwidth}
            \centering
            \includegraphics[width=.95\linewidth]{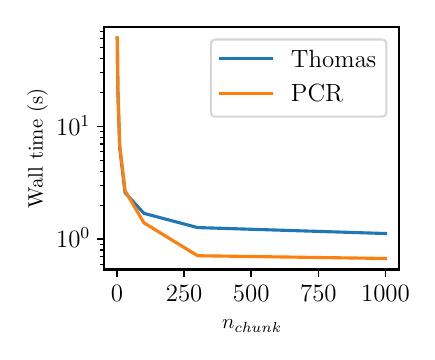}
            \caption{Wall time}
        \end{subfigure}%
        \begin{subfigure}{.375\textwidth}
            \centering
            \includegraphics[width=.95\linewidth]{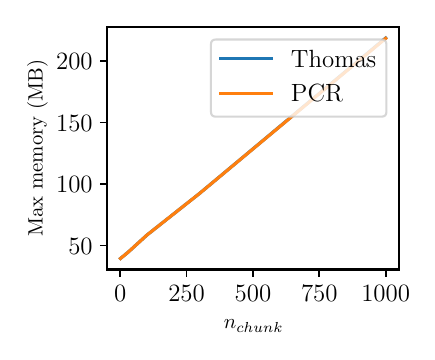}
            \caption{Memory}
        \end{subfigure}
    \end{subfigure}

    \emph{Neural ODE}
    \begin{subfigure}{1.0\textwidth}
        \centering
        \begin{subfigure}{.375\textwidth}
            \centering
            \includegraphics[width=.95\linewidth]{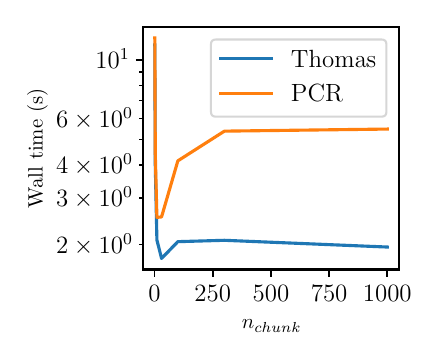}
            \caption{Wall time}
        \end{subfigure}%
        \begin{subfigure}{.375\textwidth}
            \centering
            \includegraphics[width=.95\linewidth]{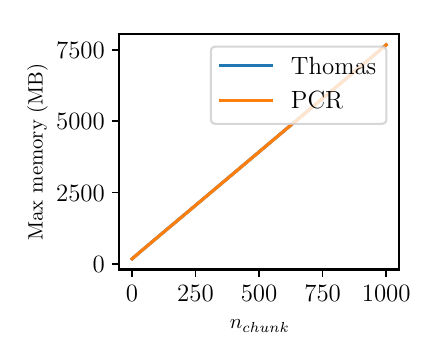}
            \caption{Memory}
        \end{subfigure}
    \end{subfigure}
    
    \caption{Chunked time integration performance study for the \emph{largest} $n_{size}$}
    \label{fig:largest}
\end{figure}

Figure \ref{fig:largest} illustrates the efficiency of the chunked time integration for all the problems, showing results for the \emph{largest} problem size $n_{size}$.  These figures plot the total wall time and the maximum memory use for each problem as a function of the time chunk size $n_{chunk}$ and for both the Thomas and PCR algorithms for solving the linearized systems.  The batched time integration algorithm achieves speedups of greater than 90x for large chunk sizes for all four problems.  It does this by trading increased memory requirements for increased bandwidth on the GPU.  As $n_{chunk}$ increases the memory required to solver the linearized systems increases, as the algorithm considers more time steps at once.  However, this also increases the amount of vectorizable work available to the GPU which decreases the time required to integrate the problem through the full $n_{time} = 2000$ time steps.  The Thomas algorithm outperforms PCR for all problems except the Chaboche model at the larger chunk sizes.

\begin{figure}
    \centering
    \emph{Mass-damper-spring}
    \begin{subfigure}{1.0\textwidth}
        \centering
        \begin{subfigure}{.375\textwidth}
            \centering
            \includegraphics[width=.95\linewidth]{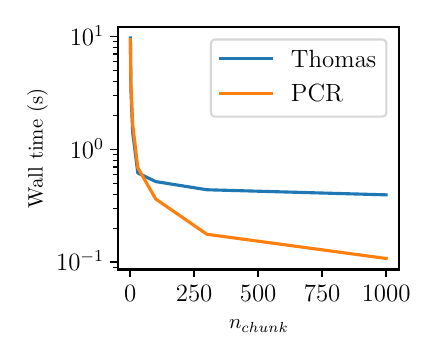}
            \caption{Wall time}
        \end{subfigure}%
        \begin{subfigure}{.375\textwidth}
            \centering
            \includegraphics[width=.95\linewidth]{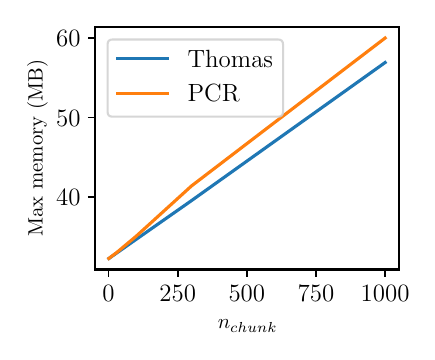}
            \caption{Memory}
        \end{subfigure}
    \end{subfigure}

    \emph{Neuron}
    \begin{subfigure}{1.0\textwidth}
        \centering
        \begin{subfigure}{.375\textwidth}
            \centering
            \includegraphics[width=.95\linewidth]{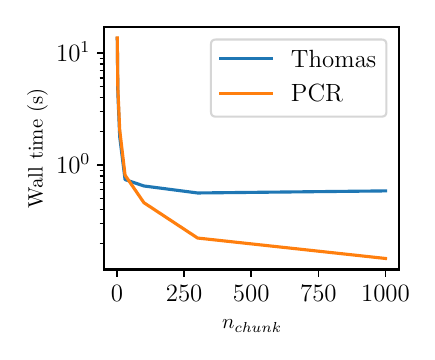}
            \caption{Wall time}
        \end{subfigure}%
        \begin{subfigure}{.375\textwidth}
            \centering
            \includegraphics[width=.95\linewidth]{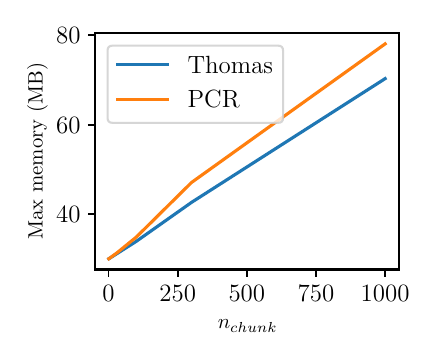}
            \caption{Memory}
        \end{subfigure}
    \end{subfigure}

    \emph{Chaboche}
    \begin{subfigure}{1.0\textwidth}
        \centering
        \begin{subfigure}{.375\textwidth}
            \centering
            \includegraphics[width=.95\linewidth]{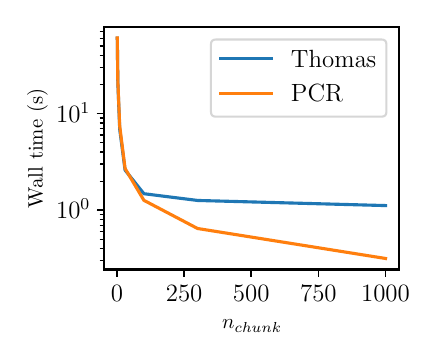}
            \caption{Wall time}
        \end{subfigure}%
        \begin{subfigure}{.375\textwidth}
            \centering
            \includegraphics[width=.95\linewidth]{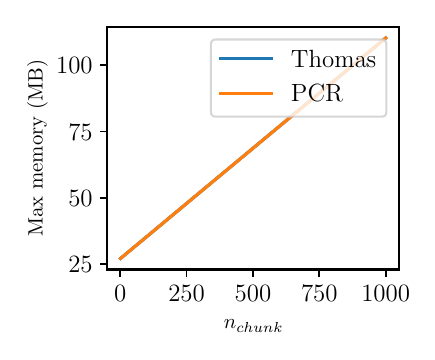}
            \caption{Memory}
        \end{subfigure}
    \end{subfigure}

    \emph{Neural ODE}
    \begin{subfigure}{1.0\textwidth}
        \centering
        \begin{subfigure}{.375\textwidth}
            \centering
            \includegraphics[width=.95\linewidth]{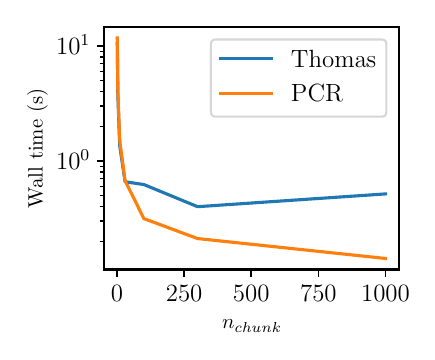}
            \caption{Wall time}
        \end{subfigure}%
        \begin{subfigure}{.375\textwidth}
            \centering
            \includegraphics[width=.95\linewidth]{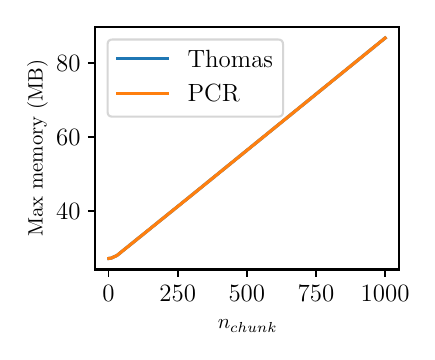}
            \caption{Memory}
        \end{subfigure}
    \end{subfigure}
    
    \caption{Chunked time integration performance study for the \emph{smallest} $n_{size}$}
    \label{fig:smallest}
\end{figure}

Figure \ref{fig:smallest} shows the same results but now for the \emph{smallest} problem size.  Now PCR outperforms Thomas for all problems for larger chunk sizes.  These results suggest Thomas's algorithm will be more effective for larger ODE systems (larger values of $n_{size}$), while PCR will be more effective for small problems.

The complexity analysis in Section \ref{sec:bidiagonal} explains this result.  Ignoring the factorization cost, which is equal for the two algorithms, the asymptotic time complexity of PCR scales with $\log_2 n_{chunk}$ of the chunk size, but with $n_{size}^3$ for the problem size.  Thomas's algorithm scales linearly with $n_{chunk}$ but quadratically with $n_{size}$.  So we expect PCR to have the advantage for small $n_{size}$, where the $\log_2$ versus linear advantage in the chunk size governs, but Thomas's algorithm to have the advantage for large $n_{size}$ where the difference between the $n_{size}^2$ scaling for Thomas's algorithm versus $n_{size}^3$ for PCR governs.

Table \ref{tab:best} shows the maximum speed up between the sequential time integration case ($n_{chunk}=1$) and the best chunk size for the vectorized time integration algorithm, for any batch size, for any problem size, and for the best linear solver algorithm, for each example problem.  All problems achieve speedup of at least 90x, with the Chaboche model reaching nearly 200x.  The performance of the chunked time integration follows the complexity of the underlying ODE, i.e. the complexity of the function $\dot{y}$ defining the ODE system.  The Chaboche model is the most complex and achieves the best speedup while the neural ODE and the mass-damper-spring system are the least complex and achieve the worst speedup.  This reflects the main advantage of vectorizing over the number of time steps in the series: the model can vectorize the evaluation of the model rate of change and Jacobian.  However, the performance of the chunked time integration algorithm is impressive for all the sample problems, even for the comparatively simple model form of the neural ODE.

\begin{table}
\centering
\begin{tabular}{ll}
\hline
Problem            & Best speedup \\ \hline
Mass-damper-spring & 106          \\
Neuron             & 94          \\
Chaboche           & 190          \\
Neural ODE         & 97           \\ \hline
\end{tabular}
\caption{Best speed up for each case between sequential and chunked time integration.}
\label{tab:best}
\end{table}

Some of the timing examples in Figs. \ref{fig:largest} and \ref{fig:smallest} demonstrate an optimal chunk size that gives the best performance.  For example, the mass-damper-spring problem with the largest problems size (Fig. \ref{fig:largest}(a)) shows an optimal chunk size of around 50 for PCR and 300 for Thomas's algorithm.  The optimal chunk size will be problem-, algorithm- (i.e. PCR versus Thomas, analytical versus AD Jacobian, etc.), data size-, and device-specific.  However, we can consider a few general factors that contribute to determining the optimal chunk size.

The total amount of parallel computation presented to the GPU depends on the product of the number of parallel data batches ($n_{batch}$) and the time chunk size ($n_{chunk}$).  There will be an optimal computation bandwidth size that depends on the device and the structure of the ODE.  For a fixed number of data, i.e. fixing $n_{batch}$, values for $n_{chunk}$ could be too small, not presenting enough parallel work to the device, or they could be too large, presenting too much work for the device to effectively handle.  This basic trade off between the $n_{chunk}$ and the amount of vectorized computational work explains the optimal chunk size for several of the examples.

As noted above, the time complexity of Thomas's algorithm versus PCR affects the optimal $n_{chunk}$ size, when comparing these two approaches.  The optimal chunk size will be smaller for PCR, compared to Thomas's algorithm, for larger problems (larger $n_{size}$) but larger for PCR compared to Thomas's algorithm for smaller problems.  The timing study data supports this trend.

Several of these timing examples do not show an optimal chunk size --- the wall time monotonically decreases as the chunk size increases.  Repeating these timing examples with longer data sets, with larger values of $n_{time}$, would allow for larger values of $n_{chunk}$, and would eventually find an optimal value of chunk size.  The wall time required to run these analyses on a single GPU for the one-at-a-time, serial time integration prevented us from finding the optimal value for all combinations of parameters considered here.

\subsection{AD versus analytic Jacobians \label{subsec:adjac}}

The main disadvantage of implicit time integration, compared to explicit methods, is the need to calculate the Jacobian matrix required to solve the discrete implicit update formula with Newton's method.  This study compares the efficiency of calculating the Jacobian analytically, with batched backward mode AD, and with batched forward mode AD.  The parameter study is the same as in the previous section, summarized in Table \ref{tab:big}, but now expanding the study to consider all three methods for calculating the Jacobian.

\begin{figure}
    \centering
    \emph{Mass-damper-spring}
    \begin{subfigure}{1.0\textwidth}
        \centering
        \begin{subfigure}{.375\textwidth}
            \centering
            \includegraphics[width=.95\linewidth]{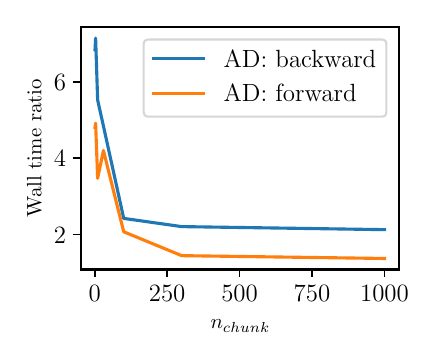}
            \caption{Wall time}
        \end{subfigure}%
        \begin{subfigure}{.375\textwidth}
            \centering
            \includegraphics[width=.95\linewidth]{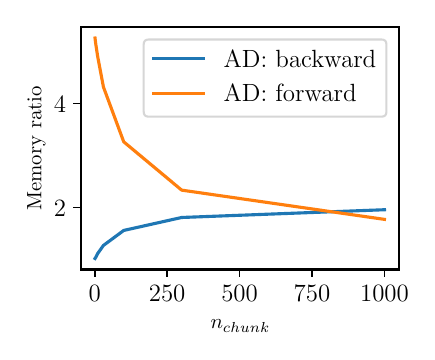}
            \caption{Memory}
        \end{subfigure}
    \end{subfigure}

    \emph{Neuron}
    \begin{subfigure}{1.0\textwidth}
        \centering
        \begin{subfigure}{.375\textwidth}
            \centering
            \includegraphics[width=.95\linewidth]{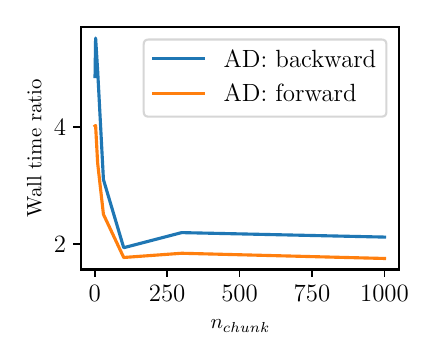}
            \caption{Wall time}
        \end{subfigure}%
        \begin{subfigure}{.375\textwidth}
            \centering
            \includegraphics[width=.95\linewidth]{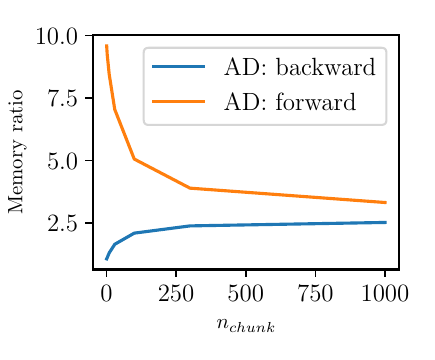}
            \caption{Memory}
        \end{subfigure}
    \end{subfigure}

    \emph{Chaboche}
    \begin{subfigure}{1.0\textwidth}
        \centering
        \begin{subfigure}{.375\textwidth}
            \centering
            \includegraphics[width=.95\linewidth]{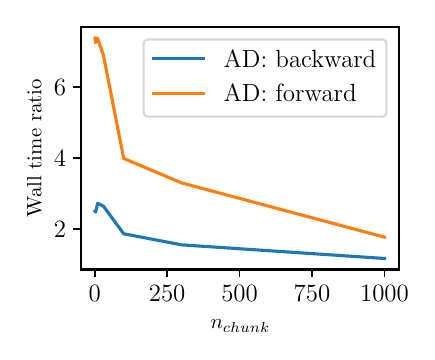}
            \caption{Wall time}
        \end{subfigure}%
        \begin{subfigure}{.375\textwidth}
            \centering
            \includegraphics[width=.95\linewidth]{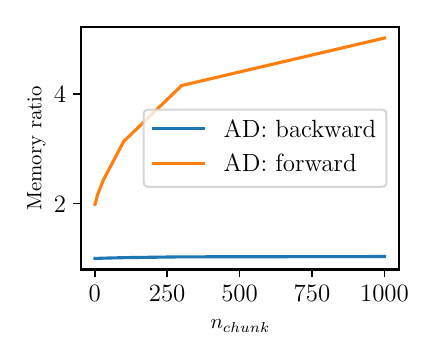}
            \caption{Memory}
        \end{subfigure}
    \end{subfigure}

    \emph{Neural ODE}
    \begin{subfigure}{1.0\textwidth}
        \centering
        \begin{subfigure}{.375\textwidth}
            \centering
            \includegraphics[width=.95\linewidth]{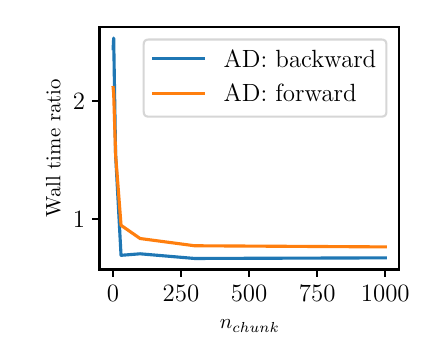}
            \caption{Wall time}
        \end{subfigure}%
        \begin{subfigure}{.375\textwidth}
            \centering
            \includegraphics[width=.95\linewidth]{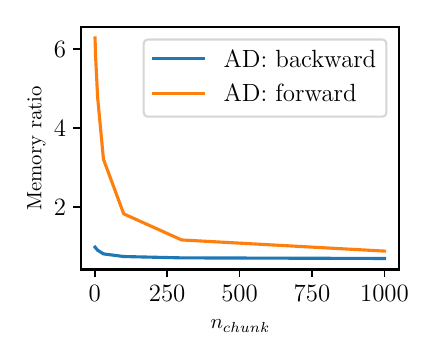}
            \caption{Memory}
        \end{subfigure}
    \end{subfigure}
    
    \caption{Cost of calculating the Jacobian with AD, compared to an analytic Jacobian, using the largest $n_{size}$ and the Thomas linear solver.}
    \label{fig:jacobian-type}
\end{figure}

Figure \ref{fig:jacobian-type} plots the ratio for the total wall time and the maximum memory use between the forward and backward mode AD Jacobian and the analytic Jacobian as a function of chunk size $n_{chunk}$.  These plots then show how much more expensive the AD calculation is compared to the analytic Jacobian for both computational time and memory use.  This figure plots results for the largest value of $n_{size}$ and for the Thomas linear solver.

Calculating the Jacobian with AD slows down the forward/backward calculation by a factor of 2 to 8, depending on the problem type and the chunk size.  Increasing the chunk size decreases the wall time penalty, compared to the analytic Jacobian, as the vectorized AD becomes comparatively more efficient.

Using AD increases the maximum required memory by a factor of 2 to 10 compared to an analytic calculation.  Again, the comparative memory penalty decreases as the chunk size increases.

Forward and backward mode AD are competitive in terms of wall time, with forward mode outperforming backward mode for the mass-damper-spring and neuron problems and backward mode winning for the Chaboche and neural ODE problems.  The timing difference between the two algorithms decreases as the chunk size increases.  AD actually outperforms the hand-coded Jacobian for the neural ODE problem for some chunk sizes.  This reflects the relatively unoptimized hard-coded Jacobian compared to the very efficient torch backward pass for this example.

The memory cost of the forward mode calculation is generally more than for backward mode AD.  The trend in terms of the problem size is ODE-specific.  For most of the examples here, the comparative disadvantage of forward mode lessens as the problem size increases.  For the mass-damper-spring system forward mode actually becomes less memory intensive for very large problems.

In general, and coding the analytic Jacobian of the ODE provides the best performance.  Failing that, batched AD provides the required Jacobian for implicit time integration with only, roughly, a 2x penalty in wall time and memory use for larger problems.  The choice of forward or backward mode AD will vary with the specific problem, though if memory is the primary concern then backward mode AD generally will be preferred.

\section{Conclusions \label{sec:conclusions}}

This paper summarizes an efficient algorithm for integrating ODEs with implicit time integration in a vectorized manner on GPUs.  The approach applies the adjoint method to efficiently calculate the derivative of the integrated time series with respect to the ODE parameters and vectorized, chunked time integration to increase the GPU throughput and improve computational efficiency for data-sparse problems.  The main conclusions of the study are:
\begin{enumerate}
    \item Parallelizing the time integration scheme by chunking some number of time steps into a single discrete time integration step can greatly improve the efficiency of both the forward and backward passes.  Wall time speed up can exceed a factor of 100x, depending on the amount of data (batch size) and the particular ODE.
    \item We can efficiently extend past work on GPU-enabled explicit time integration to implicit numerical integration algorithms.  The conclusions reached with explicit integration approaches mostly hold true for implicit methods.  Specifically, the adjoint approach for calculating gradients of the integrated time series with respect to the model parameters strictly dominates AD methods, both in time and memory use.
    \item Solving the time-discrete equations for implicit time integration approaches with Newton's method requires the Jacobian of the ODE with respect to the state.  Vectorized AD, here implemented with the \texttt{functorch} library, can efficiently provide this Jacobian for cases where deriving and implementing an analytic Jacobian is difficult.  Wall time and memory penalties are only roughly a factor of two, compared to an analytic Jacobian, for realistically large problems.
    \item We compare Thomas's algorithm to parallel cyclic reduction for solving the linearized equations resulting from applying Newton's method to the backward Euler approach.  The optimal method varies with the number of coupled ODEs in the model --- PCR is more efficient for smaller models while Thomas's algorithm is more efficient for larger models.  Asymptotic complexity analyses of the two algorithms explain this numerical result.
\end{enumerate}

We provide an open source implementation of chunked time integration for both the forward and backward Euler methods in the open-source \emph{pyoptmat} python package.  This package provides an \texttt{ode.odeint\_adjoint} method mimicking the \texttt{scipy} \texttt{odeint} interface, but implemented with GPU support through \texttt{torch} and with options to control the time chunk size.

The methods describe here provide a viable means for training ODE models for problems with sparse data, even for stiff problems.  Future work might focus specifically on neural ODE models for such problems, with an emphasis on network architectures, initialization strategies, and constraints to provide efficient, physically-reasonable models.

\section*{Acknowledgements}

The research was sponsored by the U.S. Department of Energy under Contract No. DE-AC02-06CH11357 with Argonne National Laboratory, managed and operated by UChicago Argonne LLC. The US Department of Energy, Office of Nuclear Energy supported this work through the Advanced Reactor Technologies program.

This manuscript has been authored by UChicago Argonne LLC under Contract No. DE-AC02-06CH11357 with the U.S. Department of Energy. The United States Government retains and the publisher, by accepting the article for publication, acknowledges that the United States Government retains a nonexclusive, paid-up, irrevocable, world-wide license to publish or reproduce the published form of this manuscript, or allow others to do so, for United States Government purposes. The Department of Energy will provide public access to these results of federally sponsored research in accordance with the DOE Public Access Plan: (http://energy.gov/downloads/doe-public-access-plan).

\bibliographystyle{unsrtnat}
\bibliography{references.bib}

\appendix

\section{Adjoint approach for implicit time integration}

\subsection{For simple time integration \label{subsec:adjoint-simple}}

The basic idea of the adjoint method is the same for both implicit and explicit time integration: calculate the gradient of the discrete time response of the integrated ODE system by integrating the adjoint ODE backward in time from the final time of interest back to the starting time.  To calculate a consistent gradient the same scheme must be used for both the forward and backward/adjoint problem.  So solving the adjoint problem for implicit, backward Euler time integration simply requires integrating the adjoint ODE backward in time with the backward Euler method.

There are several derivations in the literature on obtaining the adjoint ODE (c.f. \cite{bradley2013pde}), which is:
\begin{equation}
    \dot{\lambda}(t^\prime) = -\lambda \frac{\partial h}{\partial y} - \frac{\partial L}{\partial y} \delta(t - t^\prime)
    \label{eq:adjoint}
\end{equation}
where we integrate this system of ODEs backwards in time from $t_{n_{time}}$ to $0$.  Once we have the time-history of the adjoint the actual loss gradient values are given by:
\begin{equation}
    g_i = \int_{t_{n_{time}}}^{t_i} \lambda \frac{\partial h}{\partial p} dt.
\end{equation}

The Dirac delta term ``jumps'' the value of the adjoint equation wherever the loss function $L$ observes the values of the integrated time series of the state, $y_i$.  In practice this means we need to integrate the adjoint equation through the same time steps for which we integrated the original, forward problem.  

Algorithm \ref{alg:scalar-adjoint} provides the discrete time integration procedure for integrating the adjoint problem and calculating the gradient at each point in the time series using backward Euler time integration.  Recall $j_i=\frac{\partial h}{\partial y_i}$ is the ODE Jacobian at step $i$.  Thus far we have not used the derivative $\frac{\partial h}{\partial p}$, which our implementation calculates with backward mode AD.

\begin{algorithm}
\caption{Backward Euler adjoint integration.}\label{alg:scalar-adjoint}
\begin{algorithmic}
\State $i \gets n_{time}$
\State $\lambda_{n_{time}} \gets 0$
\While{$i \ge 0$}
\State $\lambda_{i-1} \gets \left(I + j_{i-1} \Delta t_{i-1} \right)^{-1} \lambda_i + \frac{\partial L}{\partial y_i}$
\State $g_{i-1} \gets \left(\lambda_i - \frac{\partial L}{\partial y_i} \right) \frac{\partial h}{\partial p}$ 
\State $i \gets i - 1$
\EndWhile
\end{algorithmic}
\end{algorithm}

Even though the forward problem $h$ may be, and nearly always is, nonlinear the adjoint problem in Eq. \ref{eq:adjoint} is linear.  This means each step of implicit time integration requires only solving a single system of linear equations, rather than a full Newton-Raphson loop.  We also do not need to retain in memory the values of the discrete, integrated adjoint $\lambda_i$ after each value is consumed by the gradient calculation.

\citet{chen2018neural} describe the explicit, forward Euler implementation of the adjoint solve, which we omit here.

\subsection{For chunked time integration \label{subsec:adjoint-chunk}}

Extending the backward, adjoint pass to chunked time integration is straightforward.  We are simply integrating a new system of ODEs (Eq. \ref{eq:adjoint}) backward in time using the same chunked time, backward Euler approach in Algorithm \ref{alg:backward-euler}.  However, the Dirac delta terms are somewhat difficult to deal with in our implementation of the backward Euler integration procedure.  Instead, we elect to deal with the chunked adjoint step as a linear algebra problem.

Starting from the discrete update formula in Algorithm \ref{alg:scalar-adjoint}
\begin{equation}
    \lambda_{i-1} = \left( I + j_{i-1} \Delta t_{i-1} \right)^{-1} \lambda_i + \frac{\partial L}{\partial y_i}
\end{equation}
and replacing 
\begin{equation}
    \lambda_{i-1} = \lambda_{i-j} = \lambda_{i} + \Delta \lambda_{j},
\end{equation}
as with the forward pass, gives
\begin{equation}
    \left(I + j_{i-j-1} \Delta t_{i-j-1} \right) \Delta \lambda_{i-j-1} - \Delta \lambda_{i-j} = \frac{\partial L}{\partial y_{i-j}} + j_{i-j-1} \Delta t_{i-j-1} \left(\frac{\partial L}{\partial y_{i-j}} - \lambda_i \right).
\end{equation}
This is another batched, blocked bidiagonal linear system that can be solved with the methods described in Section \ref{sec:bidiagonal}.  For the results here, and our reference implementation in \emph{pyoptmat}, we use the same solution strategy (i.e. Thomas's algorithm or PCR) as was used for the forward pass.

The formulation for the forward Euler adjoint pass is similar.

\end{document}